\documentclass{article}
 \usepackage{amssymb,amsmath,amsfonts}
\numberwithin{equation}{section}

\newtheorem{lemma}{Lemma}[section]
 \newtheorem{theorem}{Theorem}[section]
 \newtheorem{proposition}{Proposition}[section]
 \newtheorem{corollary}{Corollary}[section]
 \newtheorem{remark}{Remark}[section]

\newcommand{\p}{\alpha}
\newcommand{\pe}{\psi_E}
\newcommand{\po}{\psi_0}
\newcommand{\R}{\mathbb{R}}
\newcommand{\C}{\mathbb{C}}
\renewcommand{\le}{\leqslant}
\renewcommand{\ge}{\geqslant}
\newcommand{\<}{\langle}
\renewcommand{\>}{\rangle}
\newcommand{\xs}{\langle x \rangle^{\sigma}}
\newcommand{\xsn}{\langle x \rangle^{-\sigma}}
\newcommand{\eht}{e^{-iH(t-\tau)}}
\newcommand{\ehs}{e^{-iH(t-s)}}
\newcommand{\ehts}{e^{-iH(\tau-s)}}
\newcommand{\xst}{\langle x \rangle^{2\sigma}}
\newcommand{\Ls}{L_{\sigma}^2}
\newcommand{\Lsn}{L_{-\sigma}^2}
\newcommand{\om}{\Omega}
\newcommand{\ipt}{(\frac{1}{2}-\frac{1}{p_2})}
\newcommand{\ipo}{(\frac{1}{2}-\frac{1}{p_1})}
\newcommand{\ip}{(\frac{1}{2}-\frac{1}{p})}

\begin{document}
\title{Asymptotic stability of ground states in 3D nonlinear Schr\" odinger equation including subcritical cases}
\author{E. Kirr\thanks{Department of Mathematics, University of Illinois at Urbana-Champaign}\ \ and \" O. M\i zrak\footnotemark[1]}
\maketitle

\begin{abstract}
 \noindent We consider a class of nonlinear Schr\"{o}dinger equation in three space
 dimensions with an attractive potential. The nonlinearity is local but rather general
 encompassing for the first time both subcritical and supercritical (in $L^2$) nonlinearities. We study the asymptotic
 stability of the nonlinear bound states, i.e. periodic in time
 localized in space solutions. Our result shows that all solutions
 with small initial data, converge to a nonlinear bound state.
 Therefore, the nonlinear bound states are asymptotically stable. The proof hinges on dispersive estimates that we
 obtain for the time dependent, Hamiltonian, linearized dynamics
 around a careful chosen one parameter family of bound states that ``shadows" the nonlinear evolution of the
 system. Due to the generality of the methods we develop we expect
 them to extend to the case of perturbations of large bound states and to other nonlinear dispersive wave type equations.
 \end{abstract}

\section{Introduction}
In this paper we study the long time behavior of solutions of the
nonlinear Schr\" odinger equation (NLS) with potential in three
space dimensions (3-d):
\begin{align}
i\partial_t u(t,x)&=[-\Delta_x+V(x)]u+g(u), \quad t\in\R,\quad x\in\R^3 \label{u} \\
u(0,x)&=u_0(x) \label{ic}
\end{align}
where the local nonlinearity is constructed from the real valued,
odd, $C^2$ function $g:\mathbb{R}\mapsto\mathbb{R}$ satisfying
\begin{equation}\label{gest}
|g''(s)|\leq C(|s|^{\alpha_1}+|s|^{\alpha_2}),\quad s\in\mathbb{R}\
0<\alpha_1\le\alpha_2<3\end{equation} which is then extended to a
complex function via the gauge symmetry:
\begin{equation}\label{gsym}
g(e^{i\theta}s)=e^{i\theta}g(s)
\end{equation}
The equation has important applications in statistical physics
describing certain limiting behavior of Bose-Einstein condensates
\cite{dgps:bec,lsy:3d,htyau:gp}.

It is well known that this nonlinear equation admits periodic in
time, localized in space solutions (bound states or solitary waves).
They can be obtained via both variational techniques
\cite{bl:i,str:sw,rw:bs} and bifurcation methods
\cite{pw:cm,rw:bs,kz:as2d}, see also next section. Moreover the set
of periodic solutions can be organized as a manifold (center
manifold). Orbital stability of solitary waves, i.e. stability
modulo the group of symmetries $u\mapsto e^{-i\theta}u,$ was first
proved in \cite{rw:bs,mw:ls}, see also \cite{gss:i,gss:ii,ss:ins}.

In this paper we show that solutions of \eqref{u}-\eqref{ic} with
small initial data asymptotically converge to the orbit of a certain
bound state, see Theorem \ref{mt}. Asymptotic stability studies of
solitary waves were initiated in the work of A. Soffer and M. I.
Weinstein \cite{sw:mc1,sw:mc2}, see also
\cite{bp:asi,bp:asii,bs:as,sc:as,gnt:as}. Center manifold analysis
was introduced in \cite{pw:cm}, see also \cite{kn:Wed}.

The main contribution of our result is to allow for subcritical and
critical ($L^2$) nonlinearities, $0<\alpha_1\le 1/3$ in
\eqref{gest}. To accomplish this we develop an innovative technique
in which linearization around a one parameter family of bound states
is used to track the solution. Previously a fixed bound state has
been used, see the papers cited in the previous paragraph. By
continuously adapting the linearization to the actual evolution of
the solution we are able to capture the correct effective potential
induced by the nonlinearity $g$ into a time dependent linear
operator. Once we have a good understanding of the semigroup of
operators generated by the time dependent linearization, see Section
\ref{se:lin}, we obtain sharper estimates for the nonlinear dynamics
via Duhamel formula and contraction principles for integral
equations, see Section \ref{se:main}. They allow us to treat a large
spectrum of nonlinearities including, for the first time, the
subcritical ones.

The main challenge is to obtain good estimates for the semigroup of
operators generated by the time dependent linearization that we use.
This is accomplished in Section \ref{se:lin}. The technique is
perturbative, and similar to the one developed by the first author
and A. Zarnescu for 2-D Schr\" odinger type operators in
\cite{kz:as2d}, see also \cite{kz:as2d2}. The main difference is
that in 3-D one needs to remove the non-integrable singularity in
time at zero of the free Schr\" odinger propagator:
$$\|e^{i\Delta t}\|_{L^{1}\mapsto L^\infty}\sim |t|^{-3/2}.$$
We do this by generalizing a Fourier multiplier type estimate first
introduced by Journ\' e, Soffer, and Sogge in \cite{kn:JSS} and by
proving certain smoothness properties of the effective potential
induced by the nonlinearity, see the Appendix.

Since our methods rely on linearization around nonlinear bound states and estimates for integral
operators we expect them to generalize to the case of large
nonlinear ground states, see for example \cite{sc:as}, or the
presence of multiple families of bound states, see for example
\cite{sw:sgs}, where it should greatly reduce the restrictions on
the nonlinearity. We are currently working on adapting the method to
other spatial dimensions. The work in 2-D is almost complete, see
\cite{kz:as2d,kz:as2d2}.

\bigskip

\noindent{\bf Notations:} $H=-\Delta+V;$

$L^p=\{f:\mathbb{R}^2\mapsto \mathbb{C}\  |\ f\ {\rm measurable\
and}\ \int_{\mathbb{R}^2}|f(x)|^pdx<\infty\},$
$\|f\|_p=\left(\int_{\mathbb{R}^2}|f(x)|^pdx\right)^{1/p}$ denotes
the standard norm in these spaces;

$<x>=(1+|x|^2)^{1/2},$ and for $\sigma\in\mathbb{R},$ $L^2_\sigma$
denotes the $L^2$ space with weight $<x>^{2\sigma},$ i.e. the space
of functions $f(x)$ such that $<x>^{\sigma}f(x)$ are square
integrable endowed with the norm
$\|f(x)\|_{L^2_\sigma}=\|<x>^{\sigma}f(x)\|_2;$

$\langle f,g\rangle =\int_{\mathbb{R}^2}\overline f(x)g(x)dx$ is the
scalar product in $L^2$ where $\overline z=$ the complex conjugate
of the complex number $z;$

$P_c$ is the projection on the continuous spectrum of $H$ in $L^2;$

$H^n$ denote the Sobolev spaces of measurable functions having all
distributional partial derivatives up to order $n$ in $L^2,
\|\cdot\|_{H^n}$ denotes the standard norm in this spaces.

\bigskip

\noindent{\bf Acknowledgements:} The authors would like to thank
Wilhelm Schlag and Dirk Hundertmark for useful discussions on this
paper. Both authors acknowledge the partial support from the NSF
grants DMS-0603722 and DMS-0707800.

\section{Preliminaries. The center manifold.}\label{se:prelim}

The center manifold is formed by the collection of periodic solutions for (\ref{u}):
\begin{equation}\label{eq:per}
  u_E(t,x)=e^{-iEt}\psi_E(x)
\end{equation}
where $E\in\mathbb{R}$ and $0\not\equiv\psi_E\in H^2(\mathbb{R}^3)$
satisfy the time independent equation:
\begin{equation}\label{eq:ev}
[-\Delta+V]\psi_E+g(\psi_E)=E\psi_E
\end{equation}
Clearly the function constantly equal to zero is a solution of (\ref{eq:ev})
but
(iii) in the following hypotheses on the potential $V$ allows for a
bifurcation
with a nontrivial, one parameter family of solutions:

\bigskip

\noindent{\bf (H1)} Assume that
\begin{itemize}
  \item[(i)] There exists $C>0$ and $\rho >3$ such that:
   \begin{enumerate}
   \item $|V(x)|\le C<x>^{-\rho},\quad {\rm for\ all}\ x\in\mathbb{R}^3;$
   \item $\nabla V\in L^p(\mathbb{R}^3)$ for some $2\le p\le\infty$
   and $|\nabla V(x)|\rightarrow 0$ as $|x|\rightarrow\infty ;$
   \item the Fourier transform of $V$ is in $L^1.$
   \end{enumerate}
  \item[(ii)] $0$ is a regular point\footnote{see
 \cite[Definition 6]{gs:dis} or $M_\mu=\{0\}$ in relation (3.1) in \cite{mm:ae}}
 of the
spectrum of
  the linear operator $H=-\Delta+V$ acting on $L^2.$
  \item [(iii)]$H$ acting on $L^2$ has exactly one
  negative eigenvalue $E_0<0$ with corresponding normalized
  eigenvector $\psi_0.$ It is well known that $\psi_0(x)$ is  exponentially decaying as
$|x|\rightarrow\infty,$ and can be
  chosen strictly positive.
\end{itemize}

\par\noindent Conditions (i)1. and (ii) guarantee the applicability of dispersive estimates of Murata \cite{mm:ae} and Goldberg-Schlag \cite{gs:dis} to the
Schr\" odinger group $e^{-iHt}.$ Condition (i)2. implies certain
regularity of the nonlinear bound states while (i)3. allow us to use
commutator type estimates, see Theorem \ref{th:gjss}. All these
 are needed to obtain estimates for the semigroup of operators generated by our time dependent linearization,
 see Theorems \ref{th:lw} and \ref{th:lp} in
 section \ref{se:lin}. In particular (i)1. implies the local well
posedness in $H^1$ of the initial value problem
(\ref{u})-(\ref{ic}), see section \ref{se:main}.

By the standard bifurcation argument in Banach spaces \cite{ln:fa}
for (\ref{eq:ev}) at $E=E_0,$ condition (iii) guarantees existence
of nontrivial solutions. Moreover, these solutions can be organized
as a $C^1$ manifold (center manifold), see \cite[section
2]{kz:as2d}. Since our main result requires, we are going to show in
what follows that the center manifold is $C^2.$ We note that for
three and higher dimensions this has been sketched in \cite{gnt:as},
however they show smoothness by formal differentiation of certain
equations without proof that at least one side has indeed
derivatives.

As in \cite{kz:as2d} we decompose the solution of (\ref{eq:ev}) in
its projection onto the discrete and continuous part of the spectrum
of $H:$
$$\psi_E=a\psi_0+h,\quad a=\langle \psi_0,\psi_E\rangle,\
h=P_c\psi_E.$$ Projecting now (\ref{eq:ev}) onto $\psi_0$ and its
orthogonal complement $={\rm Range}\ P_c$ we get:
{\setlength\arraycolsep{2pt}
\begin{eqnarray}
0&=&h+(H-E)^{-1}P_cg(a\psi_0+h)\label{eq:evc}\\
0&=&E-E_0- a^{-1}\langle\psi_0,g(a\psi_0+h)\rangle\label{eq:evp}
\end{eqnarray}}
Although we are using milder hypothesis on $V$ the argument in the Appendix
of
\cite{pw:cm} can be easily adapted to show that:
$$\mathcal{F}(E,a,h)=h+(H-E)^{-1}P_cg(a\psi_0+h)$$
is a $C^2$ function from $(-\infty,0)\times\mathbb{C}\times
\Ls\cap H^2$ to $L^2_\sigma\cap H^2$ and $\mathcal{F}(E_0,0,0)=0,$ $D_h\mathcal{F}(E_0,0,0)=I.$ Therefore the implicit
function theorem applies to equation (\ref{eq:evc}) and leads to the
existence of $\delta_1>0$ and the $C^2$ function $\tilde h(E,a)$
from $(E_0-\delta_1,E_0+\delta_1)\times \{a\in\mathbb{C}\ :\
|a|<\delta_1\} $ to $L^2_\sigma\cap H^2$ such that (\ref{eq:evc})
has a unique solution $h=\tilde h(E,a)$ for all $E\in
(E_0-\delta_1,E_0+\delta_1)$ and $|a|<\delta_1.$ Note that, by gauge
invariance, if $(a,h)$ solves (\ref{eq:evc}) then
$(e^{i\theta}a,e^{i\theta}h),\ \theta\in[0,2\pi )$ is also a
solution, hence by uniqueness we have:
\begin{equation}\label{eq:hsym}
\tilde h(E,a)=\frac{a}{|a|}\tilde h(E,|a|).
\end{equation}
Because $\psi_0$ is real valued, we could apply the implicit function theorem
to
(\ref{eq:evc}) under the restriction $a\in\mathbb{R}$ and $h$ in the subspace
of
real valued functions as it is actually done in \cite{pw:cm}. By uniqueness of
the
solution we deduce that $\tilde h(E,|a|)$ is a real valued function.

Consider now the restriction of $\tilde h(E,a)$ to $a\in\mathbb{R},\
|a|<\delta_1.$ This is now a real valued $C^2$ function on
$(E_0-\delta_1,E_0+\delta_1)\times (-\delta_1,\delta_1)$ which, by
\eqref{eq:hsym}, is odd in the second variable. We now differentiate
\eqref{eq:evc} with $h=\tilde h(E,a),$ to obtain the following
estimates for the first and second derivatives of $\tilde h$ on
$(E,a)\in(E_0-\delta_1,E_0+\delta_1)\times (-\delta_1,\delta_1):$
 \begin{eqnarray}
 \frac{\partial\tilde h}{\partial a}(E,a)&=&-(D_h\mathcal{F})^{-1}(E,a,\tilde
 h(E,a))[(H-E)^{-1}P_c g'(a\psi_0+\tilde h(E,a))\psi_0]=\mathcal{O}(|a|^{1+\alpha_1})\nonumber\\
\frac{\partial\tilde h}{\partial E}(E,a)&=&(D_h\mathcal{F})^{-1}(E,a,\tilde
 h(E,a))[(H-E)^{-2}P_c g(a\psi_0+\tilde h(E,a))]=\mathcal{O}(|a|^{2+\alpha_1})\nonumber\\
 \frac{\partial^2\tilde h}{\partial a^2}(E,a)&=&
 -(D_h\mathcal{F})^{-1}\left[(H-E)^{-1}P_c g''(a\psi_0+\tilde h(E,a))\left(\psi_0+\frac{\partial\tilde h}{\partial a}\right)^2\right]=\mathcal{O}(|a|^{\alpha_1}) \nonumber \end{eqnarray}

\begin{eqnarray}
 \frac{\partial^2\tilde h}{\partial E\partial a}(E,a)=\mathcal{O}(|a|^{1+\alpha_1})&=&
 (D_h\mathcal{F})^{-1}\left[(H-E)^{-2}P_c g'(a\psi_0+\tilde h)\left(\psi_0+\frac{\partial\tilde h}{\partial a}\right)\right]\nonumber\\
 &-&(D_h\mathcal{F})^{-1}\left[(H-E)^{-1}P_c g''(a\psi_0+\tilde h)\left(\psi_0+\frac{\partial\tilde h}{\partial a}\right)\frac{\partial\tilde h}{\partial E}\right]\nonumber\\
  \frac{\partial^2\tilde h}{\partial E^2}(E,a)=\mathcal{O}(|a|^{2+\alpha_1})&=&
  +(D_h\mathcal{F})^{-1}\left[(H-E)^{-1}P_c g''(a\psi_0+\tilde h)\left(\frac{\partial\tilde h}{\partial E}\right)^2\right]\nonumber\\
  &-&(D_h\mathcal{F})^{-1}\left[2(H-E)^{-3}P_c g(a\psi_0+\tilde h)-2(H-E)^{-2}P_c g'(a\psi_0+\tilde h)\frac{\partial\tilde h}{\partial E}\right]\nonumber
 \end{eqnarray}
where we used $D_h\mathcal{F}(E,a,\tilde
 h(E,a))$ is invertible with bounded inverse and $D_h\mathcal{F}(E,0,0)=I,$ $(H-E)^{-1}$ is bounded
and analytic operator in $E\in(E_0-\delta_1,E_0+\delta_1),$ and
$g'(s)=\mathcal{O}(s^{1+\alpha_1}),$ $g''(s)=\mathcal{O}(s^{\alpha_1})$ as
$s\rightarrow 0.$

Replacing now $h=\tilde h(E,a),\
(E,a)\in(E_0-\delta_1,E_0+\delta_1)\times (-\delta_1,\delta_1) $ in
(\ref{eq:evp}) we get:
\begin{equation}\label{eq:evp1}
E-E_0=a^{-1}\langle\psi_0,g(a\psi_0+\tilde h(E,a))\rangle.
\end{equation}
To this we can apply again the implicit function theorem by
observing that $ G(E,a)=E-E_0-a^{-1}\langle\psi_0,g(a\psi_0+\tilde
h(E,a))\rangle$ is a $C^1$ function from
$(E_0-\delta_1,E_0+\delta_1)\times (-\delta_1,\delta_1)$ to
$\mathbb{R}$ with the properties $G(E_0,0)=0,$ $\partial_E
G(E_0,0)=1.$ We obtain the existence of $0<\delta\le\delta_1,$ and
the $C^1$ even function $\tilde E:(-\delta , \delta)\mapsto
(E_0-\delta,E_0+\delta)$  such that, for $|E-E_0|,|a|<\delta,$ the
unique solution of (\ref{eq:evp}) with $h=\tilde h(E,a),$ is given
by the $E=\tilde E(a).$ Note that $\tilde E$ is $C^2$ except at
$a=0$ because $G$ is $C^2$ except at $a=0,$ and:
\begin{eqnarray}
\frac{d\tilde E}{da}(a)&=&-\frac{\partial_a G(E(a),a)}{\partial_E
G(E(a),a)}=\mathcal{O}(|a|^{\alpha_1})\nonumber\\
\frac{d^2\tilde E}{da^2}(a)&=&\mathcal{O}(|a|^{\alpha_1-1})\qquad {\rm
for}\ a\not=0,\ {\rm recall\ that\ }0<\alpha_1\le 1.\nonumber
\end{eqnarray}

If we now define the odd function:
$$h(a)\equiv\tilde h(E(a),a),\qquad -\delta< a <\delta$$
we get a $C^2$ function because, for $a\not=0,$ based on the
previous estimates on the derivatives of $\tilde h$ and $\tilde E,$
we have
 $$
\frac{d^2h}{da^2}(a)=\frac{\partial\tilde h}{\partial
E}\frac{d^2\tilde E}{da^2}+\frac{\partial^2\tilde h}{\partial
E^2}\left(\frac{d\tilde E}{da}\right)^2+2\frac{\partial^2\tilde
h}{\partial E\partial a}\frac{d\tilde E}{da}+\frac{\partial^2\tilde
h}{\partial a^2} =\mathcal{O}(|a|^{\alpha_1}),
 $$
 hence, by L'Hospital
 $$\frac{d^2h}{da}(0)\stackrel{def}{=}\lim_{a\rightarrow 0}\frac{\frac{dh}{da}(a)-0}{a}=\lim_{a\rightarrow
 0}\frac{d^2h}{da^2}(a)=0.$$

We now extend $h$ to complex values via the rotational symmetry
\eqref{eq:hsym}:
 $$
 h(a)=\frac{a}{|a|}\tilde h(E(|a|),|a|).
 $$ We have just proved:

\begin{proposition}\label{pr:cm} There exist $\delta>0,$
the $C^2$ function
$$h:\{a\in\mathbb{R}\times\mathbb{R}\ :\ |a|<\delta\}\mapsto L^2_\sigma\cap H^2
,$$ and the $C^1$ function $E:(-\delta,\delta)\mapsto\mathbb{R}$
such that for $|E-E_0|<\delta$ and $|\langle\psi_0,\psi_E\rangle
|<\delta$ the eigenvalue problem (\ref{eq:ev}) has a unique solution
up to multiplication with $e^{i\theta},\ \theta\in [0,2\pi),$ which
can be represented as a center manifold:
 \begin{equation}\label{eq:cm}
 \psi_E=a\psi_0+h(a),\ E=E(|a|), \quad \langle\psi_0,h(a)\rangle =0,\quad h(e^{i\theta}a)=e^{i\theta}h(a),
|a|<\delta .\end{equation} Moreover
$E(|a|)=\mathcal{O}(|a|^{1+\alpha_1})$,
$h(a)=\mathcal{O}(|a|^{2+\alpha_1}),$ and for $a\in\mathbb{R},\
|a|<\delta,$ $h(a)$ is a real valued function with
$\frac{d^2h}{da^2}(a)=\mathcal{O}(|a|^{\alpha_1})$
\end{proposition}

Since $\psi_0(x)$ is exponentially decaying as $|x|\rightarrow\infty$ the
proposition implies that $\psi_E\in L^2_\sigma .$ A regularity argument, see
\cite{sw:mc1}, gives a stronger result:

\begin{corollary}\label{co:decay} For any $\sigma\in\mathbb{R},$
there exists a finite constant $C_\sigma$ such that:
$$\|<x>^\sigma\psi_E\|_{H^2}\le C_\sigma\|\psi_E\|_{H^2}.$$
\end{corollary}

\begin{remark}\label{rmk:reg} By standard regularity methods, see for example \cite[Theorem 8.1.1]{caz:bk}, one can show $\psi_E\in H^3.$
Hence by Sobolev imbeddings both $\psi_E$ and $\nabla\psi_E$ are continuous and converge to zero as $|x|\rightarrow\infty.$
\end{remark}

\begin{remark}\label{rmk:pos} By standard variational methods, see for example \cite{rw:bs}, one can show that the real valued solutions of
\eqref{eq:ev} do not change sign. Then Harnack inequality for
$H^2\bigcap C(\mathbb{R}^3)$ solutions of \eqref{eq:ev} implies that
these real solution cannot take the zero value. Hence $\psi_E$ given
by \eqref{eq:cm} for $a\in\mathbb{R}$ is either strictly positive or
strictly negative.
\end{remark}

In section \ref{se:lin} we also need some smoothness for the
effective (linear) potential induced by the nonlinearity which
modulo rotations of the complex plane is given by:
$$Dg|_{\pe}[u+iv]=g'(\pe)u+i\frac{g(\pe)}{\pe}v,\qquad \psi_E\ge 0$$
namely:

\bigskip

\noindent{\bf (H2)} Assume that for the positive solution of
\eqref{eq:ev} we have $\widehat{g'(\pe)},\
\widehat{\frac{g(\pe)}{\pe}}\in L^1(\mathbb{R}^3)$ where $\hat{f}$
stands for the Fourier transform of the function $f.$
\bigskip

\noindent In concrete cases the hypothesis may be checked directly
using the regularity of $\psi_E,$ the solution of an uniform
elliptic e-value problem. In general we can prove the following
result:

\begin{proposition}\label{pr:dg} If the following holds

\bigskip

\noindent{\bf (H2')} $g$ restricted to reals
 has third derivative except at zero and $|g'''(s)|<\frac{C}{s^{1-\alpha_1}}+Cs^{\alpha_2 -1},$ $s>0,$  $0<\alpha_1\le\alpha_2;$

\bigskip

\noindent then for the nonnegative solution of \eqref{eq:ev}, $\pe,$ we have
$\widehat{g'(\pe)}\in L^1$ and $\widehat{\frac{g(\pe)}{\pe}}\in
L^1$.
\end{proposition}

We will give the proof in the Appendix.

We are going to decompose the solution of \eqref{u}-\eqref{ic} into
a projection onto the center manifold and a correction. For orbital
stability the projection which minimizes the $H^1$ norm of the
correction is used, see for example \cite{mw:ls}, while for
asymptotic stability one wants to remove periodic in time components
of the correction. Currently there are two different ways to
accomplish this. First and most used one is to keep the correction
orthogonal to the discrete spectrum of a fixed linear Schr\" odinger
operator ``close" to the dynamics, see \cite{kz:as2d,pw:cm}. For
example in \cite{kz:as2d} the linear Schr\" odinger operator is
$-\Delta + V$ and the correction is always orthogonal on its sole
eigenvector $\psi_0,$ hence the decomposition becomes
$$u=a\psi_0+h(a)+correction,\qquad {\rm where}\ a=\langle
\psi_0,u\rangle.$$ Second technique is to use the invariant
subspaces of the actual linearized dynamics at the projection, see
for example \cite{gnt:as}. While more complicated the latter is the
only one capable to render our main result. Since there are slight
mistakes in the previous presentations of this decomposition we are
going to describe it in what follows.

Consider the linearization of \eqref{u} at function on the
center manifold $\psi_E=a\psi_0+h(a),\ a=a_1+ia_2\in\mathbb{C},\
|a|<\delta:$
 \begin{equation}\label{eq:ldE}
 \frac{\partial w}{\partial t}=-iL_{\psi_E}[w]-iEw
 \end{equation}
where
 \begin{equation}\label{def:linop}
 L_{\psi_E}[w]=(-\Delta+V-E)w+Dg_{\psi_E}[w]=(-\Delta+V-E)w+\lim_{\varepsilon\in\mathbb{R},\ \varepsilon\rightarrow
 0}\frac{g(\psi_E+\varepsilon w)-g(\psi_E)}{\varepsilon}
 \end{equation}

{\bf Properties of the linearized operator}:
 \begin{enumerate}
 \item $L_{\psi_E}$ is real linear and symmetric with respect to the
 real scalar product $\Re\langle\cdot,\cdot\rangle,$ on
 $L^2(\mathbb{R}^3),$ with domain $H^2(\mathbb{R}^3).$
 \item Zero is an e-value for $-iL_{\psi_E}$ and its generalized
 eigenspace includes $\left\{\frac{\partial\psi_E}{\partial a_1},\frac{\partial\psi_E}{\partial
 a_2}\right\}$
 \end{enumerate}

The real linearity of $L_{\psi_E}$ follows from \eqref{def:linop}.
For symmetry consider first the case of a real valued
$\psi_E=a\psi_0+h(a),\ a\in (-\delta,\delta)\subset\mathbb{R}.$ Then
 for $w=u+iv\in H^2(\mathbb{R}^3),\ u,v$ real valued we have
$$L_{\psi_E}[u+iv]=L_+[u]+iL_-[v]$$ with $L_+[u],$ $L_-[v]$ being
real valued and symmetric:
 \begin{eqnarray}
 L_+[u]&=&(-\Delta+V-E)u+g'(\psi_E)u\nonumber\\
 L_-[v]&=&(-\Delta+V-E)v+\frac{g(\psi_E)}{\psi_E}v.\nonumber
 \end{eqnarray}
To determine the expression for $L_-$ we used the rotational
symmetry \eqref{gsym}:
 $$
 g(e^{i\theta}\psi_E)=e^{i\theta}g(\psi_E)
 $$
and we differentiate it with respect to $\theta$ at $\theta=0$ to
get
 \begin{equation}\label{symgen}
 Dg_{\psi_E}[i\psi_E]=ig(\psi_E).
 \end{equation}
Now,
 $$\Re\langle L_{\psi_E}[u+iv],u_1+iv_1\rangle=\Re\langle L_+[u],u_1\rangle
+\Re\langle L_-[v],v_1\rangle=\Re\langle u,L_+[u_1]\rangle
+\Re\langle v,L_-[v_1]\rangle=\Re\langle
u+iv,L_{\psi_E}[u_1+iv_1]\rangle$$ hence $L_{\psi_E}$ is symmetric
for real valued $\psi_E.$

For a complex valued function on the center manifold
$\psi_E=a\psi_0+h(a),\ a\in\mathbb{C},\ |a|<\delta$ there exists
$\theta\in [0,2\pi)$ such that $a=|a|e^{i\theta}$ and
$$\psi_E=e^{i\theta}(|a|\psi_0+h(|a|))=e^{i\theta}\psi_E^{real}$$
where $\psi_E^{real}$ is real valued and on the center manifold.
Using again the rotational symmetry of $g$ \eqref{gsym} we get:
 \begin{equation}\label{simlinop}
 L_{\psi_E}[w]=e^{i\theta}L_{\psi_E^{real}}[e^{-i\theta}w].
 \end{equation}
Since $e^{i\theta}$ is a unitary linear operator on the real Hilbert
space $L^2(\mathbb{R}^3)$ and, due to the argument above,
$L_{\psi_E^{real}}$ is symmetric we get that $L_{\psi_E}$ is
symmetric.

For the second property, we observe that substituting $w=i\psi_E$ in
\eqref{def:linop} and using \eqref{symgen}, \eqref{eq:ev} we get
$$L_{\psi_E}[i\psi_E]=i[(-\Delta+V-E)\psi_E+g(\psi_E)]=0.$$
Hence zero is an e-value for $-iL_{\psi_E}$ and
$i\frac{\psi_E}{|a|}$ for $a\not=0$ and $i\psi_0=\lim_{a\rightarrow
0}i\frac{\psi_E}{a}$ for $a=0$ are the corresponding eigenvectors.
Moreover by differentiating \eqref{eq:ev} with respect to $a_1=\Re
a\in\mathbb{R}$ or $a_2=\Im a\in\mathbb{R}$ we get
$$-iL_{\psi_E}\left[\frac{\partial\psi_E}{\partial a_j}\right]=-\frac{\partial E}{\partial a_j}i\psi_E,\quad j=1,2.$$
Since $\frac{\partial E}{\partial a_j}=E'(|a|)\frac{\partial
|a|}{\partial a_j}\in\mathbb{R}$ we deduce that
$\frac{\partial\psi_E}{\partial a_j},\ j=1,2$ are in the generalized
eigenspace of zero. Note that, by differentiating
$h(e^{i\theta}a)=e^{i\theta}h(a)$ with respect to $\theta$ at
$\theta=0$ we get $Dh|_a[ia]=ih(a)$ and, via \eqref{eq:cm},
$D\psi_E|_a[ia]=i\psi_E.$ Since the differential can be written with
the help of the gradient:
$$i\psi_E=D\psi_E|_a[ia]=\frac{\partial\psi_E}{\partial
 a_1}\Re [ia]+\frac{\partial\psi_E}{\partial
 a_2}\Im [ia],$$ we infer that
 $$i\psi_E\in {\rm span} \left\{\frac{\partial\psi_E}{\partial a_1},\frac{\partial\psi_E}{\partial
 a_2}\right\}\quad {\rm or\ equivalently}\ \psi_E\in {\rm span} \left\{i\frac{\partial\psi_E}{\partial a_1},i\frac{\partial\psi_E}{\partial
 a_2}\right\}$$ where the span is taking over the reals\footnote{One can actually show that, for small $|a|,$ zero
is the only e-value of $-iL_{\psi_E}$ and the corresponding
eigenspace is spanned by $\frac{\partial\psi_E}{\partial a_j},\
j=1,2.$ However this is not needed in our argument.}.

One can now decompose $L^2(\mathbb{R}^3)$ into invariant subspaces
with respect to $-iL_{\psi_E}$:

 $$L^2(\mathbb{R}^3)={\rm span} \left\{\frac{\partial\psi_E}{\partial a_1},\frac{\partial\psi_E}{\partial
 a_2}\right\}\oplus
 {\cal H}_a.$$
The standard choice is to use the projection along the dual basis:
$${\cal H}_a=\{\phi_1,\phi_2\}^\perp$$
where the orthogonality is with respect to the real scalar product,
and $\phi_1,\ \phi_2$ are in the generalized eigenspace of the
adjoint of $-iL_{\psi_E}$ corresponding to the eigenvalue zero, and
$\phi_1$ is orthogonal to $\frac{\partial\psi_E}{\partial a_2}$ but
not to $\frac{\partial\psi_E}{\partial a_1}$ while $\phi_2$ is
orthogonal to $\frac{\partial\psi_E}{\partial a_1}$ but not to
$\frac{\partial\psi_E}{\partial a_2}.$ Since $L_{\psi_E}$ is
symmetric we have $(-iL_{\psi_E})^*=L_{\psi_E}i$ and a direct
calculations shows that one can choose
 $$\phi_1=-i\frac{\partial\psi_E}{\partial a_2},\qquad
 \phi_2=i\frac{\partial\psi_E}{\partial a_1}$$
as long as $\Re\langle i\frac{\partial \psi_E}{\partial
a_1},\frac{\partial \psi_E}{\partial a_2}\rangle\not=0.$ But
 $$\Re\langle i\frac{\partial \psi_E}{\partial
 a_1},\frac{\partial \psi_E}{\partial
a_2}\rangle=\Re\langle i\psi_0,i\psi_0\rangle=1,\qquad {\rm at}\ a=0
 $$
and since $\psi_E$ is $C^2$ in $a_1,\ a_2$ we have:
\begin{remark} By possible choosing $\delta>0$ smaller than the one
in Proposition \ref{pr:cm} we get:
 \begin{equation}\label{jacobi:a0}\Re\langle i\frac{\partial \psi_E}{\partial
 a_1},\frac{\partial \psi_E}{\partial
a_2}\rangle=\Re\langle i\psi_0,i\psi_0\rangle\ge\frac{1}{2}.
 \end{equation}
\end{remark}
Consequently, for $|a|<\delta,$
 \begin{equation}\label{def:ha}
 {\cal H}_a=\left\{-i\frac{\partial\psi_E}{\partial a_2},i\frac{\partial\psi_E}{\partial
 a_1}\right\}^\perp,\qquad {\rm and}\ L^2(\mathbb{R}^3)={\rm span} \left\{\frac{\partial\psi_E}{\partial a_1},\frac{\partial\psi_E}{\partial
 a_2}\right\}\oplus
 {\cal H}_a.
 \end{equation}

Our goal is to decompose the solution of \eqref{u} at each
time into:
 $$u=\psi_E+\eta=a\psi_0+h(a)+\eta,\qquad \eta\in{\cal H}_a$$
which insures that $\eta$ is not in the non-decaying directions
(tangent space of the central manifold) ${\rm span}
\left\{\frac{\partial\psi_E}{\partial
a_1},\frac{\partial\psi_E}{\partial
 a_2}\right\}$ of the linearized equation
\eqref{eq:ldE} around $\psi_E$. The fact that this can be done in an
unique manner is a consequence of the following lemma\footnote{This
is an immediate consequence of the implicit function theorem but we
find the proof in \cite{gnt:as} to be incomplete.}:

\begin{lemma}\label{lem:decomp} There exists $\delta_1>0$ such that
 any $\phi\in L^2(\mathbb{R}^3)$ satisfying $\|\phi\|_{L^2}\le\delta_1$
can be uniquely decomposed:
 $$\phi =\psi_E+\eta=a\psi_0+h(a)+\eta$$
where $a=a_1+ia_2\in\mathbb{C},\ |a|<\delta,\ \eta\in {\cal H}_a.$
Moreover the maps $\phi\mapsto a$ and $\phi\mapsto \eta$ are $C^1$
and there exist constant $C$ independent on $\phi$  such that
$$|a|\le 2\|\phi\|_{L^2},\qquad \|\eta\|_{L^2}\le C\|\phi\|_{L^2}.$$
\end{lemma}

\smallskip\par{\bf Proof:} Consider the map
$F:\{a=(a_1,a_2)\in\mathbb{R}^2\ :\ |a|<\delta\}\times
L^2(\mathbb{R}^3)\mapsto\mathbb{R}\times\mathbb{R}:$
\begin{equation}\label{eq:jacobi} F(a_1,a_2,\phi)=\left(\Re\langle -i\frac{\partial
\psi_E}{\partial
 a_2},\phi-\psi_E\rangle,\Re\langle i\frac{\partial \psi_E}{\partial
 a_1},\phi-\psi_E\rangle\right)\end{equation}
 where $\psi_E=a\psi_0+h(a),\ a=a_1+ia_2.$ Since $h(a)$ is $C^2,$
 $F$ is a $C^1$ map and:
 \begin{eqnarray}
 F(0,0,0)&=&0\nonumber\\
 \frac{\partial F}{\partial
 (a_1,a_2)}(0,0,0)&=&\mathbb{I}_{\mathbb{R}^2}\nonumber
 \end{eqnarray}
where for the calculation of the Jacobi matrix we used
\eqref{jacobi:a0}.

The implicit function theorem implies that there exist $\delta_2\le
\delta$ and a $C^1$ map:
$$\tilde F=(\tilde F_1,\tilde F_2):B(0,\delta_2)\subset
L^2(\mathbb{R}^3)\mapsto\mathbb{R}\times\mathbb{R}$$ such that the
only solutions of
$$F(a_1,a_2,\phi)=0$$
in $|a|=|a_1+ia_2|<\delta_2,\ \|\phi\|_{L^2}<\delta_2$ are given by
$$(a_1=\tilde F_1(\phi),a_2=\tilde F_2(\phi),\phi).$$
Now, for an arbitrary $\phi\in B(0,\delta_2)\subset
L^2(\mathbb{R}^3),$ since
$$\phi =\psi_E+\eta=a\psi_0+h(a)+\eta$$
with $a=a_1+ia_2\in\mathbb{C},\ |a|<\delta_2\le\delta,\ \eta\in
{\cal H}_a$ is equivalent to $F(a_1,a_2,\phi)=0$ we get that there
is a unique choice:
$$a_1=\tilde F_1(\phi),\quad a_2=\tilde F_2(\phi),\quad
\eta=\phi-a\psi_0-h(a).$$ Moreover, by choosing
$\delta_1\le\delta_2$ such that
 $$\|D\tilde F_\phi\|\le 2\qquad \forall\phi\in
L^2(\mathbb{R}^3),\ \|\phi\|_{L^2}\le\delta_1$$ where the norm is
the operator norm from $L^2(\mathbb{R}^3)$ into
$\mathbb{R}\times\mathbb{R},$ we get, for all $\phi\in
L^2(\mathbb{R}^3),\ \|\phi\|_{L^2}\le\delta_1:$
$$|a|=\sqrt{a_1^2+a_2^2}\le 2\|\phi\|_{L^2}$$
and
$$\|\eta\|_{L^2}\le \|\phi\|_{L^2}+\|\psi_E\|_{L^2}\le
\|\phi\|_{L^2}+|a|+\|h(a)\|_{L^2}\le C\|\phi\|_{L^2}$$ where $C\ge
3+2\sup_{a\in\mathbb{C},|a|\le\delta_2}\|Dh_a\|.$ Note that the
existence of $\delta_1$ is insured by the continuity of $D\tilde F$
and, from the implicit function theorem:
$$D\tilde F_0=D_\phi F|_{\phi=0}$$
and the latter has norm one being the projection operator onto
$\psi_0.$

This finishes the proof of Lemma \ref{lem:decomp}. $\Box$

\begin{remark}\label{rmk:h-1}
Both the decomposition \eqref{def:ha} and Lemma \ref{lem:decomp} can
be extended without modifications to $H^{-1}(\mathbb{R}^3)$ the dual
of $H^1$ because $\frac{\partial \psi_E}{\partial a_j}\in H^1,\
j=1,2.$ In this case $\langle u, \phi\rangle $ denotes the
evaluation of the functional $\phi\in H^{-1} $ at $u\in H^1.$
\end{remark}

We need one more technical result relating the spaces ${\cal H}_a$
and the space corresponding to the continuous spectrum of
$-\Delta+V:$
 \begin{lemma}\label{le:pcinv}
 There exists $\delta>\delta_2>0$ such that for any $a\in\mathbb{C},\ |a|\le\delta_2$ the linear map $P_c|_{{\cal H}_a}:{\cal H}_a\mapsto {\cal
 H}_0$ is invertible, and its inverse $R_a :{\cal H}_0\mapsto {\cal
 H}_a$ satisfies:
 \begin{eqnarray}
 \|R_a\zeta\|_{L^2_{-\sigma}}&\le
 &C_{-\sigma}\|\zeta\|_{L^2_{-\sigma}},\qquad \sigma\in\mathbb{R}\ {\rm and\ for\ all}\ \zeta\in {\cal H}_0\cap L^2_{-\sigma}\label{raest1}\\
 \|R_a\zeta\|_{L^p}&\le
 &C_p\|\zeta\|_{L^p},\qquad 1\le p<\infty\ {\rm and\ for\ all}\ \zeta\in {\cal H}_0\cap
 L^p\label{raest2}\\
 \overline{R_a\zeta}&=&R_a\overline\zeta\label{racc}
 \end{eqnarray}
where the constants $C_{-\sigma},\ C_p>0$ are independent of
$a\in\mathbb{C},\ |a|\le\delta_2.$
 \end{lemma}

\smallskip\par{\bf Proof:} Since $\psi_0$ is orthogonal to ${\cal H}_0,$ by continuity we can choose
$\delta>\tilde\delta_2>0$ such that $\psi_0\notin {\cal H}_a$ for
$|a|<\tilde\delta_2.$ Consequently $P_c|_{{\cal H}_a}$ is one to
one, otherwise from $\phi\in {\cal H}_a,\ \phi\neq 0,\ P_c\phi=0$ we
get $\phi=z\psi_0$ for some $z\in\mathbb{C},\ z\neq 0$ which
contradicts $\psi_0\notin {\cal H}_a.$

Next, for $|a|<\tilde\delta_2$ we construct $R_a :{\cal H}_0\mapsto
{\cal H}_a$ such that:
 \begin{equation}\label{pcinv}
 P_cR_a\zeta=\zeta,\qquad \forall\zeta\in {\cal H}_0.
 \end{equation}
Since $P_c$ is the projection onto $\{\psi_0\}^\perp,$ condition
\eqref{pcinv} is equivalent to
 \begin{equation}\label{radef} R_a\zeta=\zeta+z\psi_0\end{equation}
for some $z\in\mathbb{C}.$ To insure that the range of $R_a$ is in
${\cal H}_a$ we impose
 \begin{equation}\label{zsys}
 \Re z\langle -i\frac{\partial\psi_E}{\partial a_2},\psi_0\rangle=-\Re\langle -i\frac{\partial\psi_E}{\partial
 a_2},\zeta\rangle,\qquad \Re z\langle i\frac{\partial\psi_E}{\partial a_1},\psi_0\rangle=-\Re\langle i\frac{\partial\psi_E}{\partial
 a_1},\zeta\rangle.
 \end{equation}
This linear system of two equations with two unknowns, $\Re z$ and
$\Im z,$ is uniquely solvable whenever $\psi_0\notin {\cal H}_a.$
Note that for $a=0$ the system becomes:
$z=\langle\psi_0,\zeta\rangle.$

In \eqref{radef} we now choose $z$ to be the unique solution of
\eqref{zsys} and obtain a well defined linear map $R_a :{\cal
H}_0\mapsto {\cal H}_a$ satisfying \eqref{pcinv}.

Consequently, $P_c|_{{\cal H}_a}$ is also onto, hence invertible and
its inverse is $R_a.$ Moreover, by the continuity of the
coefficients of \eqref{zsys} with respect to $a$ we can choose
$\delta_2\le\tilde\delta_2$ such that, for all $|a|\le\delta_2:$
 \begin{equation}\label{z:est}|z|\le 2\ \sqrt{(\Re\langle -i\frac{\partial\psi_E}{\partial
 a_2},\zeta\rangle)^2+(\Re\langle i\frac{\partial\psi_E}{\partial
 a_1},\zeta\rangle)^2}.\end{equation} Hence, via \eqref{radef} and H\" older inequality we get:
 $$\|R_a\zeta\|_Y\le\|\zeta\|_Y+2\|\psi_0\|_Y\|\zeta\|_Y\sqrt{\left\|\frac{\partial\psi_E}{\partial
 a_2}\right\|_{Y^*}^2+\left\|\frac{\partial\psi_E}{\partial
 a_1}\right\|_{Y^*}^2},$$
which, for the choice $Y=L^2_{-\sigma}(\mathbb{R}^3),\
Y^*=L^2_{-\sigma}(\mathbb{R}^3)$ respectively $Y=L^p(\mathbb{R}^3),\
Y^*=L^{p'}(\mathbb{R}^3),\ \frac{1}{p}+\frac{1}{p'}=1$ give
\eqref{raest1}, respectively \eqref{raest2}. The constants are
independent of $a$ due to the continuous dependence of
$\frac{\partial \psi_E}{\partial a_j},\ j=1,2$ on $a\in\mathbb{C}$
in the compact $|a|\le\delta_2,$ and their exponential decay in
time, see proposition \ref{pr:cm} and corollary \ref{co:decay}.

Now, $P_c$ commutes with complex conjugation because it is the
orthogonal projection onto ${\psi_0}^\perp$ and $\psi_0$ is real
valued. Then \eqref{racc} follows from $R_a$ being the inverse of
$P_c.$

The proof of Lemma \ref{le:pcinv} is now complete. $\Box$

 We are now ready to prove our main result.

\section{Main Result}\label{se:main}
\begin{theorem}\label{mt}
Assume that the nonlinear term in \eqref{u} satisfies \eqref{gest}
and \eqref{gsym}. In addition assume that hypothesis (H1) and either
(H2) or (H2') hold. Let $p_1=3+\alpha_1,\ p_2=3+\alpha_2$. Then
there exists an $\varepsilon_0$ such that for all initial conditions
$u_0(x)$ satisfying
$$\max\{\|u_0\|_{L^{p_2'}},\|u_0\|_{H^1}\}\leq\varepsilon_0,\qquad \frac{1}{p_2'}+\frac{1}{p_2}=1$$
the initial value problem (\ref{u})-(\ref{ic}) is globally
well-posed in $H^1$ and the solution decomposes into a radiative
part and a part that asymptotically converges to a ground state.

More precisely, there exist a $C^1$ function
$a:\mathbb{R}\mapsto\mathbb{C}$ such that, for all $t\in\R$ we have:
\begin{equation}
u(t,x)=\underbrace{a(t)\psi_0(x)+h(a(t))}_{\psi_E(t)}+\eta(t,x)
\label{dc} \end{equation} where $\psi_E(t)$ is on the central
manifold (i.e it is a ground state) and $\eta(t,x)\in {\cal
H}_{a(t)},$ see Proposition \ref{pr:cm} and Lemma \ref{lem:decomp}.
Moreover there exists the ground states states
$\psi_{E_{\pm\infty}}$ and the $C^1$ function $\theta:\R\mapsto \R$
such that $\lim_{|t|\rightarrow\infty}\theta(t)=0$ and:
$$\lim_{t\rightarrow\pm\infty}\|\psi_E(t)-e^{-it(E_\pm-\theta(t))}\psi_{E_{\pm\infty}}\|_{H^2\bigcap
L^2_\sigma}=0,$$ while $\eta$ satisfies the following decay
estimates:
 \begin{eqnarray}
 \|\eta(t)\|_{L^2}&\leq &C_0(\alpha_1,\alpha_2)\varepsilon_0\nonumber\\
 \|\eta(t)\|_{L^{p_1}}&\leq &
 C_1(\alpha_1,\alpha_2)\frac{\varepsilon_0}{(1+|t|)^{3(\frac{1}{2}-\frac{1}{p_1})}},\quad p_1=3+\alpha_1\nonumber
 \end{eqnarray}
and, for $p_2=3+\alpha_2:$
\begin{enumerate}
\item [(i)] if  $\alpha_1\geq \frac{1}{3}$ or $\frac{1}{3} >\alpha_1 >
\frac{2\alpha_2}{3(3+\alpha_2)}$ then
 $$\|\eta(t)\|_{L^{p_2}}\leq
 C_2(\alpha_1,\alpha_2)\frac{\varepsilon_0}{(1+|t|)^{3(\frac{1}{2}-\frac{1}{p_2})}}
 $$
\item [(ii)] if $\alpha_1 =
\frac{2\alpha_2}{3(3+\alpha_2)}$ then
$$\|\eta(t)\|_{L^{p_2}}\leq
 C_2(\alpha_1,\alpha_2)\varepsilon_0\frac{\log(2+|t|)}{(1+|t|)^{3(\frac{1}{2}-\frac{1}{p_2})}}
 $$
\item [(iii)] if $\alpha_1 <
\frac{2\alpha_2}{3(3+\alpha_2)}$ then
 $$\|\eta(t)\|_{L^{p_2}}\leq
 C_2(\alpha_1,\alpha_2)\frac{\varepsilon_0}{(1+|t|)^{\frac{1+3\alpha_1}{2}}}
 $$
\end{enumerate}

where the constants $C_0,\ C_1$ and $C_2$ are independent of
$\varepsilon_0$.
\end{theorem}
\begin{remark}
 Note that the critical and supercritical cases $\frac{1}{3}\leq\p_1< 3$ are contained in
 $(i)$. Our results for these cases are stronger than the ones in
 \cite{pw:cm,sw:mc1,sw:mc2} because we do not require the initial
 condition to be in $L^2_\sigma,\ \sigma>1.$ Compared to
 \cite{gnt:as} we have sharper estimates for the asymptotic decay to the ground state but we
 require the initial data to be in $L^{p_2'}.$ To the best of our knowledge the
 subcritical case
 $\alpha_1<1/3$ has not been treated previously.
\end{remark}
\begin{remark}
 One can obtain estimates for the radiative part $\eta$ in $L^p$, $2\leq p\leq p_1=3+\p_1$, or $p_1\leq p\leq p_2=3+\p_2$ by Riesz-Thorin interpolation between
 $L^2$ and $L^{p_1}$ respectively between $L^{p_1}$ and $L^{p_2}.$
\end{remark}

\noindent\textbf{Proof of Theorem \ref{mt}} It is well known that
under hypothesis (H1)(i) the initial value problem (1)-(2) is
locally well posed in the energy space $H^1$ and its $L^2$ norm is
conserved, see for example \cite[Cor. 4.3.3 at p. 92]{caz:bk}.
Global well posedness follows via energy estimates from
$\|u_0\|_{H^1}$ small, see \cite[Remark 6.1.3 at p. 165]{caz:bk}.

We choose $\varepsilon_0\le \delta_1$ given by Lemma
\ref{lem:decomp}. Then, for all times, $\|u(t)\|_{L^2}\le\delta_1$
and we can decompose the solution into a solitary wave and a
dispersive component as in \eqref{dc}:
$$u(t)=a(t)\psi_0+h(a(t))+\eta(t)=\psi_E(t)+\eta(t)$$ Moreover, by possible making $\varepsilon_0$ smaller we can insure that
that $\|u(t)\|_{L^2}\le\varepsilon_0$ implies $|a(t)|\le\delta_2,\
t\in\mathbb{R}$ where $\delta_2$ is given by Lemma \ref{le:pcinv}.
In addition,
 since
  $$u\in C(\mathbb{R},H^{1}(\mathbb{R}^3))\cap
  C^1(\mathbb{R},H^{-1}(\mathbb{R}^3)),$$
and $u\mapsto a$ respectively $u\mapsto \eta$ are $C^1,$ see Remark
\ref{rmk:h-1}, we get that $a(t)$ is $C^1$ and $\eta\in
C(\mathbb{R},H^{1})\cap
  C^1(\mathbb{R},H^{-1}).$

The solution is now described by the $C^1$ function
$a:\mathbb{R}\in\C$ and $\eta(t)\in C(\R,H^1)\cap C^1(\R,H^{-1}).$
To obtain their equations we plug in (\ref{dc}) into (\ref{u}). Then
we get
\begin{align}
\frac{\partial\eta}{\partial t}+D\pe|_a a'&=-i(L_{\pe}+E)\eta-Ei\pe-iF_2(\pe,\eta) \label{eq:udecomp}
\end{align}
where $L_{\pe}$ is defined by \eqref{def:linop}
$$L_{\pe}\eta=(-\Delta+V-E)\eta-i\frac{d}{d\varepsilon}g(\pe+\varepsilon\eta)|_{\varepsilon=0}$$
and $F_2(\pe,\eta)$ denotes the nonlinear terms in $\eta$
\begin{equation}
F_2(\psi_E,\eta)=g(\psi_E+\eta)-g(\psi_E)-\underbrace{\frac{d}{d\varepsilon}g(\pe+\varepsilon\eta)|_{\varepsilon=0}}_{F_1(\psi_E,\eta)}
\label{g} \end{equation} Then projecting \eqref{eq:udecomp} onto the
invariant subspaces of $-iL_{\psi_E},$  ${\cal H}_{a},$ see
\eqref{def:ha} and the ${\rm span} \{\frac{\partial\psi_E}{\partial
a_1},\frac{\partial\psi_E}{\partial
 a_2}\}$, we obtain the equations for $\eta(t)$ and
$a(t):$
\begin{align}
 \frac{\partial\eta}{\partial t}&=-i(L_{\pe}+E)\eta-iF_2(\pe,\eta)-\tilde{F}_2(\pe,\eta) \label{eq:eta} \\
D\pe|_a a'&=-Ei\pe+\tilde{F}_2(\pe,\eta) \label{eq:a}
\end{align}
where \begin{equation}
       \tilde{F}_2(\pe,\eta)=\underbrace{\frac{\Re\<-i\frac{\partial\pe}{\partial a_2},-iF_2(\pe,\eta)\>}{\Re\<-i\frac{\partial\pe}{\partial a_2},\frac{\partial\pe}{\partial a_1}\>}}_{\beta_1(\pe,\eta)}\cdot\frac{\partial\pe}{\partial a_1}+\underbrace{\frac{\Re\<i\frac{\partial\pe}{\partial a_1},-iF_2(\pe,\eta)\>}{\Re\<i\frac{\partial\pe}{\partial a_1},\frac{\partial\pe}{\partial a_2}\>}}_{\beta_2(\pe,\eta)}\cdot\frac{\partial\pe}{\partial a_2} \label{eq:f2t}
      \end{equation}

\noindent In order to obtain the estimates for $\eta(t)$, we analyze \eqref{eq:eta}. The linear part of \eqref{eq:eta} is:
\begin{align}
\frac{\partial\zeta}{\partial t}&=-i(L_{\pe}+E)\zeta=(-\Delta+V)\zeta-i\frac{d}{d\varepsilon}g(\pe(t)+\varepsilon\zeta)|_{\varepsilon=0} \label{zeta} \\
\zeta(s)&=v \nonumber
\end{align}
Define $\Omega(t,s)v=\zeta(t)$. Then using Duhamel's principle \eqref{eq:eta} becomes
\begin{equation}
\eta(t)=\om(t,0)\eta(0)-\int_0^t \om(t,s)[iF_2(\pe,\eta)+\tilde{F}_2(\pe,\eta)]ds \label{eq:duheta}
\end{equation}
It is here where we differ from the approach \cite{sc:as,pw:cm,sw:mc1,sw:mc2}. The right-hand side of our equation contains only nonlinear terms in $\eta$. However the challenge is to obtain good dispersive estimates for the propagator $\om(t,s)$ of the linearization \eqref{zeta}, see Theorems \ref{th:lw} and \ref{th:lp}.

In order to apply a contraction mapping argument for
(\ref{eq:duheta}) we use the following Banach spaces. Let
$p_1=3+\p_1$ and $p_2=3+\p_2$,  $$Y_i=\{u\in L^2\cap L^{p_1}\cap
L^{p_2}:\sup_t
(1+|t|)^{3(\frac{1}{2}-\frac{1}{p_1})}\|u\|_{L^{p_1}}<\infty, \sup_t
\frac{(1+|t|)^{n_i}}{[\log(2+|t|)]^{m_i}}\|u\|_{L^{p_2}}<\infty,
\sup_t \|u\|_{L^2}<\infty \}$$ endowed with the norm $$
\|u\|_{Y_i}=\max\{\sup_t
(1+|t|)^{3(\frac{1}{2}-\frac{1}{p_1})}\|u\|_{L^{p_1}}, \sup_t
\frac{(1+|t|)^{n_i}}{[\log(2+|t|)]^{m_i}}\|u\|_{L^{p_2}}, \sup_t
\|u\|_{L^2} \}$$ for $i=1,2,3$, where
$n_1=n_2=3(\frac{1}{2}-\frac{1}{p_2})$, $n_3=\frac{1+3\p_1}{2}$,
$m_1=m_3=0$ and $m_2=1$.

Consider the nonlinear operator in (\ref{eq:duheta}):
$$N(u)=\int_0^t \om(t,s)[iF_2(\pe,u)+\tilde{F}_2(\pe,u)]ds$$
\begin{lemma}\label{lm}
Consider the cases:
$$ 1.\ \alpha_1\geq \frac{1}{3}\ {\rm or}\  \frac{1}{3} >\alpha_1 >
\frac{2\alpha_2}{3(3+\alpha_2)}; \quad 2.\ \alpha_1 =
\frac{2\alpha_2}{3(3+\alpha_2)}; \quad 3.\ \alpha_1 <
\frac{2\alpha_2}{3(3+\alpha_2)}.$$ Then, for each case number i:
 $N : Y_i\rightarrow Y_i$ is well defined, and locally Lipschitz,
 i.e. there exists $\tilde{C_i}>0$, such that $$\|Nu_1
-Nu_2\|_{Y_i}\leq\tilde{C_i}(\|u_1\|_{Y_i}+\|u_2\|_{Y_i}+\|u_1\|_{Y_i}^{1+\p_1}+\|u_2\|_{Y_i}^{1+\p_1}+\|u_1\|_{Y_i}^{1+\p_2}+\|u_2\|_{Y_i}^{1+\p_2})\|u_1
-u_2\|_{Y_i}. $$
\end{lemma}
Note that the Lemma gives the estimates for $\eta $ in the Theorem
\ref{mt}. Indeed, if we denote:$$v=\Omega(t,0)\eta(0),$$ then
$$\|v\|_{Y_i}\leq C_0\|\eta(0)\|_{L^{p_2'}\cap H^1},$$ where
$C_0=\max\{C,C_p\}$, see theorem \ref{th:lw}. We choose $\epsilon_0$
in the hypotheses of theorem \ref{mt}, such that
$$C_0\epsilon_0\le\frac{1}{2}\Big(\sqrt{1+2/\tilde{C_i}}-1\Big)$$
Then by continuity there exists $0\le Lip\le 1$ such that:
$$\|v\|_{Y_i}\leq\frac{2-Lip}{4}\Big(\sqrt{1+2Lip/\tilde{C_i}}-1\Big).$$
Let $R=L\|v\|_{Y_i}/(2-Lip)$ and $B(v,R)$ be the closed ball in
$Y_i$ with center $v$ and radius $R$. A direct calculation shows
that the right-hand side of (\ref{eq:duheta}): $$Ku=v+Nu$$ leaves
$B(v,R)$ invariant, i.e. $K:B(v,R)\mapsto B(v,R)$, and it is a
contraction with Lipschitz constant $Lip$ on $B(v,R)$.

By the contraction mapping argument, (\ref{eq:duheta}) has a unique
solution in $Y_i$. We now have two solutions of (\ref{eq:eta}), one
in $C(\R,H^1)$ from classical well posedness theory and one in
$C(\R,L^2\cap L^{p_1}\cap L^{p_2})$, $p_1=3+\p_1$, $p_2=3+\p_2$ from
the above argument. Using uniqueness and the continuous embedding of
$H^1$ in $L^2\cap L^{p_1}\cap L^{p_2}$, we infer that the solutions
must coincide. Therefore, the time decaying estimates in the spaces
$Y_{1-3}$ hold also for the $H^1$ solution.

\noindent\textbf{Proof of Lemma \ref{lm}} Let $u_1,u_2$ be in one of
the spaces $Y_i,\ i=1,2,3.$ Then at each $s\in\R$ we have:
\begin{align}
F_2(\pe(s),u_1(s))-&F_2(\pe(s),u_2(s))=g(\pe+u_1)-g(\pe+u_2)-F_1(\pe,u_1)+F_1(\pe,u_2) \nonumber \\
&=\int_0^1\Big[\frac{d}{d\tau}g(\pe+u_2+\tau(u_1-u_2))-\frac{d}{d\tau}g(\pe+\tau(u_1-u_2))|_{\tau=0}\Big]d\tau \nonumber \\
&=\int_0^1\int_0^1\frac{d}{ds}\cdot\frac{d}{d\tau}g\big(\pe+s(u_2+\tau(u_1-u_2))\big)ds\; d\tau \nonumber \end{align}
Using the hypothesis \eqref{gest} we have $|g(u)|\leq C(|u|^{2+\p_1}+|u|^{2+\p_2})$, then taking the derivatives with respect to $\tau$ and $s$ and estimating the integral we get:
\begin{align}
|F_2(\pe,u_1)-F_2(\pe,u_2)|&\leq C\big[\underbrace{(|\pe|^{\p_1}+|\pe|^{\p_2})(|u_1|+|u_2|)|u_1-u_2|}_{A_1}\label{def:a13} \\
&+\underbrace{(|u_1|^{1+\p_1}+|u_2|^{1+\p_1})|u_1-u_2|}_{A_2}+\underbrace{(|u_1|^{1+\p_2}+|u_2|^{1+\p_2})|u_1-u_2|}_{A_3}\big].\nonumber
\end{align} By \eqref{eq:f2t} and H\" older inequality, for any $1\le q\le \infty$ we have:
 \begin{eqnarray} \|\tilde{F}_2(\pe,u_1)-\tilde{F}_2(\pe,u_2)\|_{L^q}&\leq &\tilde
C \Big(\Big\|\frac{\partial\pe}{\partial
a_2}\Big\|_{L^{p_2}}\Big\|\frac{\partial\pe}{\partial
a_1}\Big\|_{L^q}+\Big\|\frac{\partial\pe}{\partial
a_1}\Big\|_{L^{p_2}}\Big\|\frac{\partial\pe}{\partial
a_2}\Big\|_{L^q}\Big)(\|A_1\|_{L^{p'_2}}+\|A_3\|_{L^{p'_2}})\nonumber \\
&+&\tilde C \Big(\Big\|\frac{\partial\pe}{\partial
a_2}\Big\|_{L^{p_1}}\Big\|\frac{\partial\pe}{\partial
a_1}\Big\|_{L^q}+\Big\|\frac{\partial\pe}{\partial
a_1}\Big\|_{L^{p_1}}\Big\|\frac{\partial\pe}{\partial
a_2}\Big\|_{L^q}\Big)\|A_2\|_{L^{p'_1}}\label{est:tildef}\\
&\leq
&C(\|A_1\|_{L^{p'_2}}+\|A_2\|_{L^{p'_1}}+\|A_3\|_{L^{p'_2}}),\nonumber
 \end{eqnarray} where the
uniform bounds on $\frac{\partial\pe}{\partial a_j}\in H^2(\R^3),\
j=1,2,$ follow from their continuous dependence on scalar $a,$ and
$|a(t)|\le\delta_2,\ t\in\R.$

Now let us consider the difference $Nu_1-Nu_2$

\begin{equation}\label{eq:nu12}
(Nu_1-Nu_2)(t)=\int_0^t \Omega(t,s) \big[iF_2(\pe
(s),u_1(s))-iF_2(\pe (s),u_2(s))+\tilde{F}_2(\pe
(s),u_1(s))-\tilde{F}_2(\pe (s),u_2(s))]ds
\end{equation}

\begin{itemize}
\item \textbf{$L^{p_2}$ Estimate :}
 \begin{align}
    \|Nu_1-Nu_2\|_{L^{p_2}}&\leq\int_0^t \|\om(t,s)\|_{L^{p'_2}\rightarrow L^{p_2}}C\Big(2\|A_1\|_{L^{p'_2}}+\|A_2\|_{L^{p'_2}}+\|A_2\|_{L^{p'_1}}+2\|A_3\|_{L^{p'_2}}\Big)ds \nonumber \end{align}
To estimate the term containing $A_1$, observe that
$$\|(|\pe|^{\p_1}+|\pe|^{\p_2})(|u_1|+|u_2|)|u_1-u_2|\|_{L^{p'_2}}\leq\||\pe|^{\p_1}+|\pe|^{\p_2}\|_{L^{\beta}}(\|u_1\|_{L^{p_2}}+\|u_2\|_{L^{p_2}})\|u_1-u_2\|_{L^{p_2}}$$
with $\frac{1}{\beta}+\frac{2}{p_2}=\frac{1}{p'_2}$. Using Theorem
\ref{th:lp} (see also Remark \ref{rmk:lp}), we have for each case
number i: and $u_1,\ u_2\in Y_i:$

\begin{align}
\int_0^t & \|\om(t,s)\|_{L^{p'_2}\rightarrow L^{p_2}}\|A_1\|_{L^{p'_2}}ds \nonumber \\
    &\leq\int_0^t \frac{C(p_2)}{|t-s|^{3\ipt}}\||\pe|^{\p_1}+|\pe|^{\p_2}\|_{L^{\beta}}\frac{[\log(2+|s|)]^{2m_i}}{(1+|s|)^{2n_i}}(\|u_1\|_{Y_i}+\|u_2\|_{Y_i})\|u_1-u_2\|_{Y_i} ds \nonumber \\
    &\leq\frac{C(p_2)C_1 C_2}{(1+|t|)^{3\ipt}}(\|u_1\|_{Y_i}+\|u_2\|_{Y_i})\|u_1-u_2\|_{Y_i} \nonumber \end{align}
where $C_2=\sup_t \frac{(1+|t|)^{n_i}}{[\log(2+|t|)]^{m_i}}\int_0^t
\frac{[\log(2+|s|)]^{2m_i}ds}{|t-s|^{3\ipt}(1+|s|)^{2n_i}}<\infty$
since $2n_i >1$ and $C_1=\sup_t
\||\pe|^{\p_1}+|\pe|^{\p_2}\|_{L^{\beta}}.$ The uniform bounds in
$t\in\R$ for  $\|\pe\|_{L^{\alpha_j\beta}}^{\alpha_j},\ j=1,2$
follow from the continuous dependence of $\pe=a(t)\po+h(a(t))\in
H^2(\R^3)$ on $a(t)$ and $|a(t)|\le \delta_2,\ t\in\R.$

To estimate the terms containing $A_2$, observe that
$$\|(|u_1|^{1+\p_1}+|u_2|^{1+\p_1})|u_1-u_2|\|_{L^{p'_1}}\leq\big(\|u_1\|_{L^{p_1}}^{1+\p_1}+\|u_2\|_{L^{p_1}}^{1+\p_1}\big)\|u_1-u_2\|_{L^{p_1}}$$
since $\frac{1}{p'_1}=\frac{2+\p_1}{p_1}$ and
 \begin{eqnarray} \|(|u_1|^{1+\p_1}+|u_2|^{1+\p_1})|u_1-u_2|\|_{L^{p'_2}}&\leq &\left(\|u_1\|_{L^{p_1}}^{\theta(1+\p_1)}\|u_1\|_{L^{2}}^{(1-\theta)(1+\p_1)}+\|u_2\|_{L^{p_1}}^{\theta(1+\p_1)}\|u_2\|_{L^{2}}^{(1-\theta)(1+\p_1)}\right) \nonumber \\
                                                     &&\times\|u_1-u_2\|^\theta_{L^{p_1}}\|u_1-u_2\|^{1-\theta}_{L^{2}}\nonumber
                                                    \end{eqnarray}
 where $\frac{1}{p'_2}=(2+\p_1)(\frac{1-\theta}{2}+\frac{\theta}{p_1}),\ 0\le\theta\le 1. $ Again using Theorem \ref{th:lp} (see also Remark \ref{rmk:lp}), we have
\begin{align}
    \int_0^t & \|\om(t,s)\|_{L^{p'_2}\rightarrow L^{p_2}}\|A_2\|_{L^{p'_2}}ds \nonumber \\
    &\leq\int_0^t \frac{C(p_2)}{|t-s|^{3\ipt}}\cdot\frac{(\|u_1\|_{Y_i}^{1+\p_1}+\|u_2\|_{Y_i}^{1+\p_1})\|u_1-u_2\|_{Y_i}}{(1+|s|)^{3(\frac{\p_1}{2}+\frac{1}{p_2})}} ds \nonumber \\
&\leq\frac{C(p_2)C_3[\log(2+|t|)]^{m_i}}{(1+|t|)^{n_i}}(\|u_1\|_{Y_i}^{1+\p_1}+\|u_2\|_{Y_i}^{1+\p_1})\|u_1-u_2\|_{Y_i}
\nonumber \end{align} where the different decay rates $n_i$ depend
on the case number in the hypotheses of this Lemma:
 \begin{enumerate}
 \item corresponds to
$3(\frac{\p_1}{2}+\frac{1}{p_2})>1$, and  $C_3=\sup_t
(1+|t|)^{3\ipt}\int_0^t
\frac{ds}{|t-s|^{3\ipt}(1+|s|)^{3(\frac{\p_1}{2}+\frac{1}{p_2})}}<\infty;$
 \item corresponds to $3(\frac{\p_1}{2}+\frac{1}{p_2})=1$, and $C_3=\sup_t
\frac{(1+|t|)^{3\ipt}}{\log(2+|t|)}\int_0^t
\frac{ds}{|t-s|^{3\ipt}(1+|s|)}<\infty;$
 \item corresponds to $3(\frac{\p_1}{2}+\frac{1}{p_2})<1$, and $C_3=\sup_t
(1+|t|)^{\frac{1+3\p_1}{2}}\int_0^t
\frac{ds}{|t-s|^{3\ipt}(1+|s|)^{3(\frac{\p_1}{2}+\frac{1}{p_2})}}<\infty.$
 \end{enumerate}

 To estimate the term containing $A_3$, observe that
 $$\|(|u_1|^{1+\p_2}+|u_2|^{1+\p_2})|u_1-u_2|\|_{L^{p'_2}}\leq\big(\|u_1\|_{L^{p_2}}^{1+\p_2}+\|u_2\|_{L^{p_2}}^{1+\p_2}\big)\|u_1-u_2\|_{L^{p_2}}$$ since $\frac{1}{p'_2}=\frac{2+\p_2}{p_2}$. Again using Theorem \ref{th:lp} (see also Remark \ref{rmk:lp}), we have
    \begin{align}
    \int_0^t & \|\om(t,s)\|_{L^{p'_2}\rightarrow L^{p_2}}\|A_3\|_{L^{p'_2}}ds \nonumber \\
    &\leq\int_0^t \frac{C(p_2)}{|t-s|^{3\ipt}}\cdot\frac{[\log(2+|s|)]^{(2+\p_2)m_i}}{(1+|s|)^{(2+\p_2)n_i}}(\|u_1\|_{Y_i}^{1+\p_2}+\|u_2\|_{Y_i}^{1+\p_2})\|u_1-u_2\|_{Y_i} ds \nonumber \\
&\leq\frac{C(p_2)C_4C_5[\log(2+|t|)]^{m_i}}{(1+|t|)^{n_i}}(\|u_1\|_{Y_i}^{1+\p_2}+\|u_2\|_{Y_i}^{1+\p_2})\|u_1-u_2\|_{Y_i} \nonumber \end{align}
    where $C_5=\sup_t \frac{(1+|t|)^{n_i}}{[\log(2+|t|)]^{m_i}}\int_0^t \frac{[\log(2+|s|)]^{(2+\p_2)m_i}ds}{|t-s|^{3\ipt}(1+|s|)^{(2+\p_2)n_i}}<\infty$ since $(2+\p_2)n_i >1$.
     \item \textbf{$L^{p_1}$ Estimate :}
     From \eqref{eq:nu12} we have
    \begin{eqnarray}
    \|Nu_1-Nu_2\|_{L^{p_1}}(t)&\le &\|\int_0^t \Omega(t,s)
[iF_2(\pe(s),u_1(s))-iF_2(\pe(s),u_2(s))]ds\|_{L^{p_1}}\nonumber\\
 &+&C\int_0^t\|\Omega(t,s)\|_{L^{p'_1}\mapsto
L^{p_1}}\|\tilde{F}_2(\pe(s),u_1(s))-\tilde{F}_2(\pe(s),u_2(s))\|_{L^{p'_1}}ds\nonumber
  \end{eqnarray}
For the second integral we use \eqref{est:tildef} with $q=p'_1$ and
the previous estimates on $A_1,\ A_2,\ A_3$ to obtain the required
bound. For the first integral moving the norm inside the integration
and applying $L^{p'_1}\mapsto L^{p_1}$ estimates for $\Omega(t,s)$
and \eqref{def:a13} for the nonlinear term would require the control
of $A_3$ in $L^{p'_1}.$ The latter, unfortunately, can no longer be
interpolated between $L^2$ and $L^{p_2}.$ To avoid this difficulty
we separate and treat differently the part of the nonlinearity
having an $A_3$ like behavior by decomposing $\R^3$ in two disjoints
measurable sets related to the inequality \eqref{def:a13}:
 $$
 V_1(s)=\{x\in\R^3\ |\
 |F_2(\pe(s,x),u_2(s,x))-F_2(\pe(s,x),u_1(s,x))|\le C
 A_3(s,x)\},\qquad
 V_2(s)=\R^3\setminus V_1(s)$$
On $V_2(s),$ using polar representation of complex numbers, we
further split the nonlinear term into:
 \begin{eqnarray}
 \lefteqn{iF_2(\pe(s,x),u_1(s,x))-iF_2(\pe(s,x),u_2(s,x))=e^{i\theta(s,x)}CA_3(s,x)}\nonumber\\
 &&+\underbrace{e^{i\theta(s,x)}[|iF_2(\pe(s,x),u_1(s,x))-iF_2(\pe(s,x),u_2(s,x))|-CA_3(s,x)]}_{G(s,x)}\nonumber\end{eqnarray}
where, due to inequality \eqref{def:a13}, $|G(s,x)|\le
C(A_1(s,x)+A_2(s,x))$ on $V_2(s).$ Then we have:
 \begin{eqnarray}\lefteqn{\int_0^t \Omega(t,s)
[iF_2(\pe(s),u_1(s))-iF_2(\pe(s),u_2(s))]ds=\int_0^t \Omega(t,s)
 (1-\chi (s))G(s)ds}\nonumber\\
 &&+\underbrace{\int_0^t \Omega(t,s)
[\chi(s)(iF_2(\pe(s),u_1(s))-iF_2(\pe(s),u_2(s)))+(1-\chi(s))e^{i\theta(s)}CA_3(s)]ds}_{I(t)},\nonumber\end{eqnarray}
where $\chi(s)$ is the characteristic function of $V_1(s).$ Now
 $$\|\int_0^t \Omega(t,s)
 (1-\chi (s))G(s)ds\|_{L^{p_1}}\le \int_0^t\| \Omega(t,s)\|_{L^{p'_1}\mapsto
 L^{p_1}}C(\|A_1(s)\|_{L^{p'_1}}+\|A_2(s)\|_{L^{p'_1}})ds$$
 and estimates as in the previous step for $A_1$ and $A_2$ give the required
 decay. For $I(t)$ we use interpolation:
        $$\|I(t)\|_{L^{p_1}}\leq\|I(t)\|_{L^2}^{1-\theta}\|I(t)\|_{L^{p_2}}^{\theta}\leq\|I(t)\|_{L^2}^{1-\theta}\left(\int_0^t \|\om(t,s)\|_{L^{p'_2}\mapsto L^{p_2}}\|A_3\|_{L^{p'_2}} ds\right)^{\theta}$$
where $\frac{1}{p_1}=\frac{1-\theta}{2}+\frac{\theta}{p_2}$. We know
from previous step that the above integral decays as
$(1+|t|)^{-3\ipt}$ and below we will show its $L^2$ norm will be
bounded. Therefore $$\sup_t
(1+|t|)^{3\ipo}\|I(t)\|_{L^{p_1}}<\infty$$ and the $L^{p_1}$
estimates are complete.
    \item \textbf{$L^2$ Estimate :}
To estimate $L^2$ norm we cannot use $L^2\rightarrow L^2$ estimate
for $\om(t,s)$ because that would force us to control
$L^{2(\alpha_2+2)}$ which cannot pe interpolated between $L^2$ and
$L^{p_2},\ p_2=\alpha_2+3.$ We avoid this by using the
decomposition:
    $$T(t,s)v=[P_c\om(t,s)-\ehs P_c]v \quad \text{i.e.}\quad \om(t,s)=R_{a(t)} T(t,s)+R_{a(t)} \ehs P_c$$
For $T(t,s)$ we will use $L^{p'}\rightarrow L^2$ estimates, see
Theorem \ref{th:lw}, while for $\ehs P_c$ we will use Stricharz
estimates $L^\infty_tL_x^2$. We will also use a decomposition of the
nonlinear term similar to the one for $L^{p_1}$ estimates that will
allow us to estimate in a different manner this time the terms
behaving like $A_2,$ see \eqref{def:a13}. All in all we have:
    \begin{align}
    \|Nu_1-Nu_2\|_{L^2}&\leq \int_0^t \|\om(t,s)\|_{L^2\rightarrow L^2}\|\tilde{F}_2(\pe,u_1)-\tilde{F}_2(\pe,u_2)\|_{L^2}ds
    \nonumber\\
    &+\|R_{a(t)}\|_{L^2\mapsto L^2}\int_0^t \|T(t,s)\|_{L^{p'_2}\rightarrow L^2}C(\|A_1\|_{L^{p'_2}}+\|A_3\|_{L^{p'_2}})ds \nonumber \\
    &+\|R_{a(t)}\|_{L^2\mapsto L^2}\int_0^t \|T(t,s)\|_{L^{p'_1}\rightarrow L^2}C\|A_2\|_{L^{p'_1}}ds \nonumber \\
    &+\|R_{a(t)}\|_{L^2\mapsto L^2}\|\int_0^t  \ehs P_c (A_1(s)+A_3(s) ds\|_{L^2}+\|R_{a(t)}\|_{L^2\mapsto L^2}\|\int_0^t \ehs P_c A_2 ds\|_{L^2}\nonumber
    \end{align}
For the first integral we use Theorem \ref{th:lp} part (i),\
\eqref{est:tildef} with $q=2$ and the estimates we have already
obtained for $A_1,\ A_2$ and $A_3.$ We deduce that this integral is
uniformly bounded in $t\in\R.$ Similarly we get uniform boundedness
of the second and third integral by using Theorem \ref{th:lw} part
(iv).

For the fourth integral we use Stricharz estimate:
    $$\sup_{t\in\R}\|\int_0^t \ehs P_c A_1ds\|_{L^2}\leq C_s\left[\left(\int_\R \|A_1(s)\|_{L^{p'_2}}^{\gamma_2'}ds\right)^{\frac{1}{\gamma_2'}}+\left(\int_\R \|A_3(s)\|_{L^{p'_2}}^{\gamma_2'}ds\right)^{\frac{1}{\gamma_2'}}\right]$$
where $\frac{1}{\gamma_2'}+\frac{1}{\gamma_2}=1,$ and
$\frac{2}{\gamma_2}=3\left(\frac{1}{2}-\frac{1}{p_2}\right).$ Using
again the estimates we obtained before for $A_1$ and $A_3.$ we get:
\begin{align}
\| A_1\|_{L^{\gamma_2'}_sL^{{p'_2}}}&\leq C_{11}\Big[\int_\R \frac{(\log(2+|s|))^{2m_i\gamma_2'}}{(1+|s|)^{2n_i\gamma_2'}}ds\Big]^{\frac{1}{\gamma_2'}}(\|u_1\|_{Y_i}+\|u_2\|_{Y_i})\|u_1-u_2\|_{Y_i} \nonumber \\
&\leq  C_{11} C_8 (\|u_1\|_{Y_i}+\|u_2\|_{Y_i})\|u_1-u_2\|_{Y_i}
\nonumber \end{align} where $C_8= \int_\R
\frac{(\log(2+|s|))^{2m_i\gamma_1'}}{(1+|s|)^{2n_i\gamma_1'}}ds<\infty$
since $2n_i\gamma'>1$ and:
\begin{align}
\| A_3\|_{L^{\gamma_2'}_sL^{{p'_2}}}&\leq C_{12}\Big[\int_\R \frac{(\log(2+|s|))^{(2+\p_2)m_i\gamma_3'}}{(1+|s|)^{(2+\p_2)n_i\gamma_3'}}ds\Big]^{\frac{1}{\gamma_3'}}(\|u_1\|_{Y_i}^{1+\p_2}+\|u_2\|_{Y_i}^{1+\p_2})\|u_1-u_2\|_{Y_i} \nonumber \\
&\leq  C_9
(\|u_1\|_{Y_i}^{1+\p_2}+\|u_2\|_{Y_i}^{1+\p_2})\|u_1-u_2\|_{Y_i}
\nonumber \end{align} where $C_9=\int_\R
\frac{(\log(2+|s|))^{(2+\p_2)m_i\gamma_3'}}{(1+|s|)^{(2+\p_2)n_i\gamma_2'}}ds<\infty$
since $(2+\p_2)n_1\gamma_3'>1$.

Similarly, for the fifth integral:
    $$\sup_{t\in\R}\|\int_0^t \ehs P_c A_2 ds\|_{L^2}\leq C_s\Big(\int_\R \|A_2\|_{L^{p_1'}}^{\gamma_1'}ds\Big)^{\frac{1}{\gamma_1'}}$$
where $\frac{1}{\gamma_1'}+\frac{1}{\gamma_1}=1,$ and
$\frac{2}{\gamma_1}=3\left(\frac{1}{2}-\frac{1}{p_1}\right).$
Furthermore we have
\begin{align}
\| A_3\|_{L^{\gamma_1'}_sL^{{p'_1}}}&\leq C_{13}\Big[\int_\R \frac{ds}{(1+|s|)^{3(2+\p_1)\gamma_2'\ipo}}\Big]^{\frac{1}{\gamma_2'}}(\|u_1\|_{Y_i}^{1+\p_1}+\|u_2\|_{Y_i}^{1+\p_1})\|u_1-u_2\|_{Y_i} \nonumber \\
&\leq C_{13}C_{10}
(\|u_1\|_{Y_i}^{1+\p_1}+\|u_2\|_{Y_i}^{1+\p_1})\|u_1-u_2\|_{Y_i}
\nonumber \end{align} where $C_{10}=\int_\R
\frac{ds}{(1+|s|)^{3(2+\p_1)\gamma_2'\ipo}}ds<\infty$ since
$3(2+\p_1)\gamma_2'\ipo>1$ .
\end{itemize}

The $L^2$ estimates are now complete and the proof of Lemma \ref{lm}
is finished. $\Box$

We now finish the proof of Theorem \ref{mt} by analyzing the
dynamics on the center manifold and showing it converges to a ground
state. Using the fact that
$$i\psi_E=D\psi_E|_a[ia]=\frac{\partial\psi_E}{\partial
 a_1}\Re [ia]+\frac{\partial\psi_E}{\partial
 a_2}\Im [ia]$$ equation \eqref{eq:a} becomes $$D\pe|_a(a'+iEa)=\frac{\partial\psi_E}{\partial
 a_1}\Re [a'+iEa]+\frac{\partial\psi_E}{\partial
 a_2}\Im [a'+iEa]=\tilde{F}_2(\pe,\eta)=\beta_1(\pe,\eta)\frac{\partial\psi_E}{\partial
 a_1}+\beta_2(\pe,\eta)\frac{\partial\psi_E}{\partial
 a_2}$$ Hence $$|a'+iEa|=\sqrt{\beta_1^2+\beta_2^2}=b(t)$$
and
$$
\left|[a(t)e^{i\int_0^t E(s)ds}]'\right|=b(t)
            $$
Since $b(t)=\sqrt{\beta_1^2+\beta_2^2}$, and
$$\beta_1\leq\left\|\frac{\partial\pe}{\partial a_2}\right\|_{L^{p_2}}(\|A_1\|_{L^{p'_2}}+\|A_3\|_{L^{p'_2}})+\left\|\frac{\partial\pe}{\partial a_2}\right\|_{L^{p_1}}\|A_2\|_{L^{p'_1}}\leq C(\|\eta\|^{2}_{L^{p_2}}+\|\eta\|^{2+\p_2}_{L^{p_2}}+\|\eta\|^{2+\p_1}_{L^{p_1}})$$
$$\beta_2\leq\left\|\frac{\partial\pe}{\partial a_1}\right\|_{L^{p_2}}(\|A_1\|_{L^{p'_2}}+\|A_3\|_{L^{p'_2}})+\left\|\frac{\partial\pe}{\partial a_1}\right\|_{L^{p_1}}\|A_2\|_{L^{p'_1}}\leq C(\|\eta\|^{2}_{L^{p_2}}+\|\eta\|^{2+\p_2}_{L^{p_2}}+\|\eta\|^{2+\p_1}_{L^{p_1}})$$
we get $0\le b(t)\le C(1+|t|)^{1+\delta}$ for some $\delta>0,$ in
each of the cases $(i)$, $(ii)$ and $(iii)$ in the Theorem \ref{mt}.
Then, for any $\varepsilon>0$ we have
\begin{equation}
\left|a(t)e^{i\int_0^{t} E(s)ds}-a(t')e^{i\int_0^{t'}
E(s)ds}\right|\le \int_{t'}^{t} b(s)ds<\varepsilon
\end{equation}
for $t,\ t'$ sufficiently large respectively sufficiently small.
Therefore $a(t)e^{i\int_0^t E(s)ds}$ has a limit when
$t\rightarrow\pm\infty$. This means
$$e^{i\int_0^t E(s)ds}\pe=a(t)e^{i\int_0^t E(s)ds}\po+e^{i\int_0^t E(s)ds}h(a(t))=a(t)e^{i\int_0^t E(s)ds}\po+h(a(t)e^{i\int_0^t E(s)ds})\rightarrow\psi_{E_{\pm\infty}}$$
Above we used $h(e^{i\theta}a)=e^{i\theta}h(a)$, see Proposition
\ref{pr:cm}. In addition $|a(t)|\rightarrow a_{\pm}$ as
$t\rightarrow\pm$ at a rate $|t|^{-\delta}.$ Since $E(s)=E(|a(s)|$
is $C^1$ in $|a|$ on $|a|\le\delta_2,$ we deduce $|E(\pm
s)-E_\pm|\le C(1+s)^{-\delta}$ for $s\ge 0$ and some constant $C>0.$
If we denote
$$\theta(\pm t)=\frac{1}{\pm t}\int_0^{\pm t}E(s)-E_\pm ds,\qquad t\ge 0$$
then $\lim_{|t|\rightarrow\infty}\theta(t)=0$ and
$$\lim_{t\rightarrow\pm}e^{it(E_\pm-\theta(t))}\psi_E(t)=\psi_{E_\pm}.$$

 This finishes the proof of Theorem \ref{mt}. $\Box$

\section{Linear Estimates}\label{se:lin}
Consider the linear Schr\" odinger equation with a potential in
three space dimensions:
\begin{align}
i\frac{\partial u}{\partial t}&=(-\Delta+V(x))u \nonumber \\
u(0)&=u_0 \nonumber \end{align}
It is known that if V satisfies hypothesis (H1) (i) and (ii) then the radiative part of the solution, i.e. its projection onto the continuous spectrum of $H=-\Delta+V$, satisfies the estimates:
\begin{equation}
    \|e^{-iHt}P_c u_0\|_{\Lsn}\leq C_M\frac{1}{|t|^{\frac{3}{2}}}\|u_0\|_{\Ls}
    \label{eq:ls}
\end{equation}
for $\sigma>1$ and some constant $C_M>0$ independent of $u_0$ and $t\in\R$, and
\begin{equation}
    \|e^{-iHt}P_c u_0\|_{L^p}\leq C_p\frac{1}{|t|^{3(\frac{1}{2}-\frac{1}{p})}}\|u_0\|_{L^{p'}}
    \label{eq:lp}
\end{equation}
for some constant $C_p>0$ depending only on $2\leq p$. The case
$p=\infty$ in (\ref{eq:lp}) is proved by Goldberg and Schlag in \cite{gs:dis}.
The conservation of the $L^2$ norm gives the $p=2$ case:
$$\|e^{-iHt}P_c u_0\|_{L^2}=\|u_0\|_{L^2}$$ The general result
(\ref{eq:lp}) follows from Riesz-Thorin interpolation.

We would like to extend these estimates to the linearized dynamics around the center manifold. We consider the linear equation, with initial data at time $s$,
\begin{align}
i\frac{d\zeta}{dt}&=H\zeta+F_1(\psi_E,\zeta) \nonumber \\
\zeta(s)&=v \nonumber \end{align}
where $F_1(\pe,\zeta)=\frac{d}{d\varepsilon}g(\pe+\varepsilon\zeta)|_{\varepsilon=0}=\frac{\partial}{\partial u}g(u)|_{u=\pe}\zeta+\frac{\partial}{\partial \bar{u}}g(u)|_{u=\pe}\overline{\zeta}$. For the sake of simpler notation, we will use $F_1(\zeta)$.

By Duhamel's principle we have: \begin{equation} \zeta(t)=\ehs v(s)-i\int_s^t \eht F_1(\zeta)d\tau \label{soz} \end{equation}

In the next theorems we will extend estimates of type (\ref{eq:ls})-(\ref{eq:lp}) to the operators $\om(t,s)$ and
$T(t,s)$ considering the fact that $\pe(t)$ is small. Recall that
 $$T(t,s)=P_c\om(t,s)-\ehs P_c\quad\textrm{i.e.}\quad \om(t,s)=R_{a(t)} T(t,s)+R_{a(t)}\ehs P_c$$
\begin{theorem}\label{th:lw}
There exists $\varepsilon_1>0$ such that for
$\|\xs\pe\|_{H^2}<\varepsilon_1$ there exist constants $C$, $C_p >0$
with the property that for any t, s $\in\R$ the followings hold:
\begin{align}
(i)\quad &\|\om(t,s)\|_{\Ls\rightarrow\Lsn}\leq\frac{C}{(1+|t-s|)^{\frac{3}{2}}} \nonumber \\
(ii) \quad &\|T(t,s)\|_{L^1\rightarrow\Lsn}\leq\left\{ \begin{array}{ll}
\frac{C}{|t-s|^{\frac{1}{2}}} & \textrm{for $s\leq t\leq s+1$} \\
\frac{C}{(1+|t-s|)^{\frac{3}{2}}} & \textrm{for $t>s+1$} \end{array} \right.
\nonumber \\
(iii)\quad &T(t,s)\in L_t^2(\R,L^2\rightarrow\Lsn)\cap L_t^\infty(\R,L^2\rightarrow\Lsn) \nonumber \\
(iv)\quad&\|\om(t,s)\|_{L^{p'}\rightarrow\Lsn}\leq\frac{C}{|t-s|^{3(\frac{1}{2}-\frac{1}{p})}}\quad\textrm{for all $2\leq p\leq L^\infty$} \nonumber \\
&\|T(t,s)\|_{L^{p'}\rightarrow\Lsn}\leq\left\{ \begin{array}{ll}
\frac{C}{|t-s|^{(\frac{1}{2}-\frac{1}{p})}} & \textrm{for $s\leq t\leq s+1$} \\
\frac{C}{(1+|t-s|)^{3(\frac{1}{2}-\frac{1}{p})}} & \textrm{for $t>s+1$} \end{array} \right.
\nonumber \end{align}

\end{theorem}

\textbf{Proof of Theorem \ref{th:lw}} Fix $s\in\R$.

$(i)$ By definition, we have $\om(t,s)v=\zeta(t)$ where $\zeta(t)$
satisfies equation \eqref{soz}. We project \eqref{soz} onto
continuous spectrum of $H=-\Delta+V:$  \begin{equation} \xi(t)=\ehs
P_c v-i\int_s^t \eht P_c F_1(R_a\xi)d\tau \label{sox} \end{equation}
where $\xi=P_c\zeta$. We are going to prove the estimate for
$P_c\Omega(t,s)$ by showing that the nonlinear equation \eqref{sox}
can be solved via contraction principle argument in an appropriate
functional space. To this extent let us consider the functional
space
$$X_1:=\{u\in C(\R,L_{-\sigma}^2(\R^3))|\sup_{t>s}(1+(t-s))^{\frac{3}{2}}\|u(t)\|_{L_{-\sigma}^2}<\infty\}$$ endowed with the norm $$\|u\|_{X_1}:=\sup_{t>s}\{(1+(t-s))^{\frac{3}{2}}\|u(t)\|_{L_{-\sigma}^2}\}<\infty$$ Note that the inhomogeneous term in (\ref{sox}) $\xi_0=\ehs P_c v$ satisfies $\xi_0\in X_1$ and
\begin{equation}
    \|\xi_0\|_{X_1}\leq C_M\|v\|_{\Ls}
    \label{eq:z0}
\end{equation}
because of (\ref{eq:ls}). We collect the $\xi$ dependent part of the right hand side of (\ref{sox}) in a linear operator $L(s):X_1\rightarrow X_1$
\begin{equation}
    [L(s)\xi](t)=-i\int_s^t \eht P_c[F_1(R_a\xi)]d\tau  \label{L} \end{equation}
We will show that $L$ is a well defined bounded operator from $X_1$ to $X_1$ whose operator norm can be made less or equal to 1/2 by choosing $\varepsilon_1$ sufficiently small. Consequently $Id-L$ is invertible and the solution of the equation (\ref{sox}) can be written as $\xi=(Id-L)^{-1}\xi_0$. In particular $$\|\xi\|_{X_1}\leq(1-\|L\|)^{-1}\|\xi_0\|_{X_1}\leq 2\|\xi_0\|_{X_1}$$
which in combination with the definition of $\om$, the definition of the norm $X_1$ and the estimate (\ref{eq:z0}), finishes the proof of $(i)$.

It remains to prove that $L$ is a well defined bounded operator from $X_1$ to $X_1$ whose operator norm can be made less than 1/2 by choosing $\varepsilon_1$ sufficiently small.
\begin{align}
\|L(s)\xi(t)\|_{L_{-\sigma}^2}\leq & \int_s^t \|\eht P_c\|_{L_{\sigma}^2\rightarrow L_{-\sigma}^2}\|F_1(R_a\xi)\|_{\Ls}d\tau \nonumber \\
\end{align}
On the other hand
$$\|F_1(R_a\xi)\|_{\Ls}\leq\|\xst(|\pe|^{1+\p_1}+|\pe|^{1+\p_2})\|_{L^{\infty}}\|R_a \xi\|_{\Lsn}\leq(\varepsilon_1^{1+\p_1}+\varepsilon_1^{1+\p_2})\|\xi\|_{\Lsn}$$ and using the last three relations, as well as the estimate (\ref{eq:ls}) and the fact that $\xi\in X_1$ we obtain that
\begin{align}
 \|L(s)\|_{X_1\rightarrow X_1}&\leq(\varepsilon_1^{1+\p_1}+\varepsilon_1^{1+\p_2}) \sup_{t>0}(1+|t-s|)^{\frac{3}{2}}\int_s^t \frac{1}{(1+|t-\tau|)^{\frac{3}{2}}} \cdot \frac{1}{(1+|\tau-s|)^{\frac{3}{2}}}{\ } d\tau \nonumber \\
&\leq(\varepsilon_1^{1+\p_1}+\varepsilon_1^{1+\p_2}) \sup_{t>0}(1+|t-s|)^{\frac{3}{2}}\frac{1}{(1+|\frac{t-s}{2}|)^{\frac{3}{2}}}\leq C(\varepsilon_1^{1+\p_1}+\varepsilon_1^{1+\p_2}) \nonumber \end{align}
Now choosing $\varepsilon_1$ small enough we get
$$\|L\|_{X_1\rightarrow X_1}<\frac{1}{2}$$ Therefore
$$\|P_c\Omega(t,s)\|_{\Ls\rightarrow\Lsn}\leq\frac{\tilde C}{(1+|t-s|)^{\frac{3}{2}}}$$
and
$$\|\Omega(t,s)\|_{\Ls\rightarrow\Lsn}\leq \|R_{a(t)}\|_{\Lsn\rightarrow\Lsn}\|P_c\Omega(t,s)\|_{\Ls\rightarrow\Lsn}\le
\frac{C}{(1+|t-s|)^{\frac{3}{2}}}$$ by Lemma \ref{le:pcinv}.

$(ii)$ Recall that \begin{equation}
    P_c\om(t,s)v=T(t,s)v+\ehs P_c v
    \label{omt}
\end{equation}
Denote: \begin{equation}
    T(t,s)v=W(t)
    \label{w}
\end{equation}
then, by plugging in (\ref{sox}), $W(t)$ satisfies the following equation:
\begin{equation}
W(t)=\underbrace{-i\int_s^t \eht P_c [F_1(R_a\ehts P_c v)]d\tau}_{f(t)}+[L(s)W](t) \label{sow} \end{equation} By definition of $T(t,s)$ (\ref{w}) it is sufficient to prove that the solution of (\ref{sow}) satisfies $$\|W(t)\|_{\Lsn}\leq\left\{ \begin{array}{ll}
\frac{C\|v\|_{L^1}}{|t-s|^{\frac{1}{2}}} & \textrm{for $s\leq t\leq s+1$} \\
\frac{C\|v\|_{L^1}}{(1+|t-s|)^{\frac{3}{2}}} & \textrm{for $t>s+1$} \end{array} \right.$$

Let us also observe that it suffices to prove this estimate only for the forcing terms $f(t)$ because then we will be able to do the contraction principle in the functional space in which $f(t)$ will be, and thus obtain the same decay for $W$ as for $f(t)$.

This time we will consider the functional space $$X_2=\{u\in C(\R,\Lsn{\R^3})| \sup_{|t-s|>1} (1+|t-s|)^\frac{3}{2}\|u\|_{\Lsn}<\infty, \sup_{|t-s|\leq1} |t-s|^\frac{1}{2}\|u\|_{\Lsn}<\infty \}$$
endowed with the norm $$\|u\|_{X_2}=\left\{ \begin{array}{ll}
\sup_t |t-s|^\frac{1}{2}\|u\|_{\Lsn} & \textrm{for $|t-s|\leq 1$} \\
\sup_t (1+|t-s|)^\frac{3}{2}\|u\|_{\Lsn} & \textrm{for $|t-s|>1$} \end{array} \right.$$
Now we will estimate $f(t)$. First we will investigate the short time behavior of this term. If $s\leq t\leq s+1$. Recall that $F_1(u)=\frac{d}{d\tau}g(\pe+\tau u)=\frac{\partial}{\partial u}g(u)|_{u=\pe}u+\frac{\partial}{\partial \bar{u}}g(u)|_{u=\pe}\overline{u}=g_{u}u+g_{\bar{u}}\overline{u}.$
 $$\|f(t)\|_{\Lsn}\leq\|\xsn\int_s^t \eht P_c F_1(R_a\ehts P_c v) d\tau\|_{L^2}$$
For the term $g_u R_a\ehts P_c v$ we have
\begin{align}
 &\|\xsn\int_s^t \|\eht P_c g_u R_a\ehts P_c v d\tau\|_{L^2}  \nonumber\\
\leq&\|\xsn\|_{L^2}\int_s^t \|\ehs e^{iH(\tau-s)} P_c g_u R_a\ehts P_c v\|_{L^\infty} d\tau \nonumber\\
&\leq \int_s^t\frac{C}{|t-s|^{\frac{3}{2}}}\sup \|\widehat{g_u}\|_{L^1}\|v\|_{L^1}d\tau\leq C\frac{\|v\|_{L^1}}{|t-s|^{\frac{1}{2}}}\sup\|\widehat{g_u}\|_{L^1} <\infty \nonumber
\end{align} and for the term $g_{\bar{u}}\overline{R}_a e^{iH(\tau-s)} P_c v$ we have
\begin{align}
 &\|\xsn\int_s^t \|\eht P_c g_{\bar{u}}\overline{R}_a e^{iH(\tau-s)} P_c v d\tau\|_{L^2}  \nonumber\\
&\leq\|\xsn\|_{L^2}\int_s^{s+\frac{t-s}{4}} \|e^{-iH(t+s-2\tau)} \ehts P_c g_{\bar{u}}\overline{R}_a e^{iH(\tau-s)} P_c v\|_{L^\infty} d\tau \nonumber\\
&\quad+\int_{s+\frac{t-s}{4}}^t \|\eht P_c\|_{\Ls\rightarrow\Lsn}\|\xs g_{\bar{u}}\overline{R}_a \ehts P_c v\|_{L^2} d\tau \nonumber \\
&\leq\int_s^{s+\frac{t-s}{4}}\frac{C}{|t+s-2\tau|^{\frac{3}{2}}}\sup \|\widehat{g_{\bar{u}}}\|_{L^1}\|v\|_{L^1}d\tau+\int_{s+\frac{t-s}{4}}^t \frac{C}{(1+|t-\tau)^{3/2}}\|\xs g_{\bar{u}}\|_{L^2} \|\ehts P_c v\|_{L^\infty} d\tau \nonumber \\
&\leq C\frac{\|v\|_{L^1}}{|t-s|^{\frac{1}{2}}}\sup (\|\widehat{g_{\bar{u}}}\|_{L^1}+\|\xs g_{\bar{u}}\|_{L^2})<\infty \nonumber  \end{align} There we used J-S-S type estimate; see Appendix for $t=\tau-s$; $|\tau-s|\leq 1$. For the long time behavior of $f(t)$, we will split this integral into three parts to be estimated differently. For $t>s+1$, $$f(t)=\underbrace{\int_s^{s+\frac{1}{2}}\cdots}_{I_{1}}+\underbrace{\int_{s+\frac{1}{2}}^{\frac{t+s}{2}}\cdots}_{I_{2}}+\underbrace{\int_{\frac{t+s}{2}}^t \cdots}_{I_{3}}$$
Then for $t>s+1$,
\begin{align} \|I_1\|_{\Lsn}&\leq\|\xsn\int_s^{s+\frac{1}{2}} \ehs P_c F_1(R_a \ehts P_c v) d\tau\|_{L^2} \nonumber \\
&\leq\|\xsn\|_{L^2}\int_s^{s+\frac{1}{2}} \|\ehs e^{iH(\tau-s)} P_c g_u R_a\ehts P_c v\|_{L^\infty} d\tau \nonumber \\
&\quad+\|\xsn\|_{L^2}\int_s^{s+\frac{1}{2}} \|e^{-iH(t+s-2\tau)} \ehts P_c g_{\bar{u}}\overline{R}_a e^{iH(\tau-s)} P_c v\|_{L^\infty} d\tau \nonumber \\
&\leq\|\xsn\|_{L^2}\int_s^{s+\frac{1}{2}} \frac{C}{|t-s|^{\frac{3}{2}}}\|e^{iH(\tau-s)} P_c g_u R_a\ehts P_c v\|_{L^1} d\tau \nonumber \\
&\quad+\|\xsn\|_{L^2}\int_s^{s+\frac{1}{2}} \frac{C}{|t+s-2\tau|^{\frac{3}{2}}}\|\ehts P_c g_{\bar{u}}\overline{R}_a e^{iH(\tau-s)} P_c v\|_{L^1} d\tau \nonumber \\
&\leq C\|\xsn\|_{L^2}\Big(\frac{1}{|t-s|^{\frac{3}{2}}}+\frac{1}{|t-s-1|^{\frac{3}{2}}}\Big)\int_s^{s+\frac{1}{2}} \sup (\|\widehat{g_u}\|_{L^1}+\|\widehat{g_{\bar{u}}}\|_{L^1})\|v\|_{L^1}d\tau \nonumber \\
&\leq C\|\xsn\|_{L^2}\sup(\|\widehat{g_u}\|_{L^1}+\|\widehat{g_{\bar{u}}}\|_{L^1})\frac{1}{(1+|t-s|)^{\frac{3}{2}}}\|v\|_{L^1} \nonumber \end{align}
For the second integral we have \begin{align}
\|I_2\|_{\Lsn}&\leq\int_{s+\frac{1}{2}}^{\frac{t+s}{2}} \|\eht P_c\|_{\Ls\rightarrow\Lsn}\|F_1(R_a\ehts P_c v)\|_{\Ls} d\tau \nonumber \\
&\leq\int_{s+\frac{1}{2}}^{\frac{t+s}{2}} \frac{C}{(1+|t-\tau|)^{\frac{3}{2}}}\|\xs|\pe|^{1+\p}\|_{L^2}\|\ehts P_c v\|_{L^\infty}d\tau \nonumber \\
&\leq\frac{C\|v\|_{L^1}}{(1+|\frac{t-s}{2}|)^{\frac{3}{2}}}\int_{s+\frac{1}{2}}^{\frac{t+s}{2}} \frac{d\tau}{|\tau-s|^{\frac{3}{2}}}\leq\frac{C\|v\|_{L^1}}{(1+|t-s|)^{\frac{3}{2}}} \nonumber \end{align}
$I_3$ is estimated similiar to $I_2$.

$(iii)$ From \eqref{sow} we have \begin{align} \xsn W(t)=&\int_s^t \xsn\eht P_c[F_1(R_a\ehts P_c v)]d\tau+\int_s^t \xsn\eht P_c[F_1(R_a W(\tau))]d\tau \nonumber \end{align}
Then
\begin{align} \|\xsn W(t)\|_{L_t^2 L_x^2}&\leq\Big\|\int_s^t \frac{C}{(1+|t-\tau|)^{3/2}} \|\xs F_1(R_a\ehts P_c v)\|_{L^2}d\tau\Big\|_{L^2_t} \nonumber \\
&+\Big\|\int_s^t \frac{C}{(1+|t-\tau|)^{3/2}}(\|\xst g_u\|_{L^\infty}+\|\xst g_{\bar{u}}\|_{L^\infty})\|\xsn W(\tau)\|_{L^2_x}d\tau \Big\|_{L^2_t} \nonumber \\
&\leq C\|K\|_{L^1}\|v\|_{L^2}+\varepsilon_1 C\|K\|_{L^1}\|\xsn W\|_{L^2_x L^2_t} \nonumber\end{align}
Where $K(t)=(1+|t|)^{-3/2}$. For the term $\xs F_1(R_a e^{-iHt} P_c v)=\xs (g_u R_a e^{-iHt} P_c v+g_{\bar{u}}\overline{R}_a e^{iHt} P_c v)$ we used $\|\xst g_u\|_{L^\infty}$ and $\|\xst g_{\bar{u}}\|_{L^\infty}$ is uniformly bounded in $t$ since $|g_u|=|g_{\bar{u}}|\leq C(|\pe|^{1+\p_1}+|\pe|^{1+\p_2})$ and the Kato smoothing estimate $\|\xsn e^{-iHt} P_c v\|_{L^2_t(\R,L^2_x)}\leq C\|v\|_{L^2_x}$. Choosing $\varepsilon_1$ small enough we get $\|\xsn W\|_{L^2_x L^2_t}<\infty$. In other words $T(t,s)\in L_t^2(\R, L^2\rightarrow\Lsn)$. And similarly \begin{align} \|\xsn W(t)\|_{L_x^2}&\leq\int_s^t \frac{C}{(1+|t-\tau|)^{3/2}} \|\xs F_1(R_a\ehts P_c v)\|_{L^2}d\tau \nonumber \\
&+\int_s^t \frac{C}{(1+|t-\tau|)^{3/2}}(\|\xst g_u\|_{L^\infty}+\|\xst g_{\bar{u}}\|_{L^\infty})\|\xsn W(\tau)\|_{L^2_x}d\tau \nonumber \\
&\leq C\|v\|_{L^2}+\varepsilon_1 C\|\xsn W\|_{L^2_x} \nonumber\end{align} This finishes the proof of $(iii)$, $T(t,s)\in L_t^2(\R, L^2\rightarrow\Lsn)\cap L_t^\infty(\R, L^2\rightarrow\Lsn).$

$(iv)$ By Riesz-Thorin interpolation between $(ii)$ and $(iii)$ (the $L_t^\infty$ part) we get the desired estimates. $\Box$

The next step is to obtain estimates for $\om(t,s)$ and $T(t,s)$ in unweighted $L^p$ spaces.
\begin{theorem}{\label{th:lp}}
Assume that $\|\xs\pe\|_{H^2} < \varepsilon_1$ (where $\varepsilon_1$ is the one used in Theorem \ref{th:lw}).
Then there exist constants $C_2$, $C_2'$ and $C_\infty$ for all $t$, $s$ $\in$ $\R$ the following estimates hold:
\begin{align}
(i)\quad &\|\om(t,s)\|_{L^2\rightarrow L^2}\leq C_2,\quad \|T(t,s)\|_{L^2\rightarrow L^2}\leq C_2 \nonumber \\
(ii)\quad &\|\om(t,s)\|_{L^1\rightarrow L^\infty}\leq\frac{C_\infty}{|t-s|^{\frac{3}{2}}}, \quad \|T(t,s)\|_{L^1\rightarrow L^\infty}\leq\left\{ \begin{array}{ll}
C_\infty|t-s|^{\frac{1}{2}} & \textrm{for $|t-s|\leq 1$} \\
\frac{C_\infty}{|t-s|^{\frac{3}{2}}} & \textrm{for $|t-s|>1$} \end{array} \right. \nonumber \\
(iii)\quad &\|T(t,s)\|_{L^{p'}\rightarrow L^2}\leq C_2',
\quad\textrm{for }\quad p=6 \nonumber\end{align}
\end{theorem}
\begin{remark}\label{rmk:lp}
 By Riesz-Thorin interpolation from $(i)$ and $(ii)$, and from $(i)$ and $(iii)$ we get
\begin{eqnarray}
\|\om(t,s)\|_{L^{p'}\rightarrow L^p}&\leq&\frac{C_p}{|t-s|^{3\ip}},\quad\textrm{for all}\quad 2\leq p\leq\infty \nonumber \\
\|T(t,s)\|_{L^{p'}\rightarrow L^p}&\leq&\left\{ \begin{array}{ll}
C_p|t-s|^{\ip} & \textrm{for $|t-s|\leq 1$} \\
\frac{C}{|t-s|^{3\ip}} & \textrm{for $|t-s|>1$} \end{array} \right.,\quad\textrm{for all}\quad 2\leq p\leq\infty \nonumber \\
\|T(t,s)\|_{L^{p'}\rightarrow L^2}&\leq &C_p,
\quad\textrm{for all}\quad 2\leq p\leq6 \nonumber
\end{eqnarray}

\end{remark}

\noindent\textbf{Proof of Theorem \ref{th:lp}} Because of the estimate (\ref{eq:lp}) and relation $P_c\om=T+\ehs P_c$, It suffices to prove the theorem for $T(t,s)$.

$(i)$ To estimate the $L^2$ norm we will use duality argument to make use of cancelations. \begin{align}
&\|f(t)\|_{L^2}^2=\< f(t),f(t)\> \nonumber \\
&=\int_s^t\int_s^t \< \eht P_c F_1(R_a\ehts P_c v),e^{-iH(t-\tau')} P_c F_1(R_a e^{-iH(\tau' -s)} P_c v)\> d\tau' d\tau  \nonumber \\
&=\int_s^t \int_s^t \< F_1(R_a\ehts P_c v),e^{-iH(\tau-\tau')} P_c F_1(R_a e^{-iH(\tau'-s)} P_c v)\> d\tau' d\tau  \nonumber \\
&=\int_s^t \int_s^t \< \xs F_1(R_a\ehts P_c v), \xsn e^{-iH(\tau-\tau')} P_c F_1(R_ae^{-iH(\tau'-s)} P_c v) \> d\tau' d\tau  \nonumber \\
&\leq\int_s^t \int_s^t \|F_1(R_a\ehts P_c v)\|_{\Ls}\|e^{-iH(\tau-\tau')} P_c F_1(R_a e^{-iH(\tau'-s)} P_c v)\|_{\Lsn} d\tau' d\tau \nonumber \\
&\leq\int_s^t \|\xs F_1(R_a\ehts P_c v)\|_{L^2}\int_s^t\frac{C}{(1+|\tau-\tau'|)^{3/2}}\|\xs F_1(R_a e^{-iH(\tau'-s)} P_c v)\|_{L^2} d\tau' d\tau\nonumber \\
&\leq C\|\xs F_1(R_a\ehts P_c v)\|_{L_{\tau}^2 L_x^2}
\Big\|\int_s^t\frac{C}{(1+|\tau-\tau'|)^{3/2}} \|\xs F_1(R_a e^{-iH(\tau'-s)} P_c v)\|_{L_x^2} d\tau\Big\|_{L^2_{\tau}}\nonumber \\
&\leq C\|K\|_{L^1}\|\xs F_1(R_a e^{-iHt} P_c v)\|_{L_t^2 L^2_x}^2\leq C\|v\|^2_{L^2}<\infty \nonumber \end{align}
At the last line, $K(t)=(1+|t|)^{-3/2}$ and  we used convolution estimate. For the term $\xs F_1(R_a e^{-iHt} P_c v)=\xs (g_u R_a e^{-iHt} P_c v+g_{\bar{u}}\overline{R}_a e^{iHt} P_c v)$ we used the Kato smoothing estimate $\|\xsn e^{-iHt} P_c v\|_{L^2_t(\R,L^2_x)}\leq C\|v\|_{L^2_x}$.
We will estimate $L^2$ norm of $L$ similiar to $f$. 
\begin{align}
&\|L(s)W\|_{L^2}^2=\<L(s)W,L(s)W\> \nonumber \\
&=\<\int_s^t \eht P_c F_1(W(\tau))d\tau,\int_s^t e^{-iH(t-\tau')} P_c F_1(W(\tau'))d\tau'\> \nonumber \\
&=\int_s^t \int_s^t \< F_1(W(\tau)),e^{-iH(\tau-\tau')} P_c F_1(W(\tau'))\>] d\tau' d\tau  \nonumber \\
&\leq\int_s^t (\|\xs g_u\|_{L^\infty}+\|\xs g_{\bar{u}}\|_{L^\infty})\|\xsn W\|_{L^2} \nonumber \\
&\quad\quad\quad\quad\quad\quad\times\int_s^t C K(\tau-\tau')(\|\xs g_u\|_{L^\infty}+\|\xs g_{\bar{u}}\|_{L^\infty})\|\xsn W\|_{L^2} d\tau' d\tau\nonumber \\
&\leq C\|\xsn W\|_{L_x^2 L_\tau^2}
\Big\|\int_s^t CK(\tau-\tau') \|\xsn W\|_{L_x^2} d\tau'\Big\|_{L^2_{\tau}}\nonumber \\
&\leq C\|K\|_{L^1}\|\xsn W\|_{L_\tau^2 L^2_x}<\infty \nonumber \end{align} By Theorem \ref{th:lw} $(iii)$, $\|\xsn W\|_{L_\tau^2 L^2_x}<\infty$.

\noindent Therefore we conclude $\|T(s,t)\|_{L^2\rightarrow L^2}\leq C$ and $\|\om(s,t)\|_{L^2\rightarrow L^2}\leq C$

$(ii)$ Let us first investigate the short time behavior of the forcing term $f(t)$. We will assume $s\leq t\leq s+1$,

\begin{align}
\|f(t)\|_{L^\infty}&=\|\int_s^t \eht P_c F_1(R_a \ehts P_c v)d\tau\|_{L^\infty} \nonumber \\
&\leq\int_s^t \|e^{-iH(t-s)}P_c\|_{L^1\rightarrow L^\infty}\|e^{iH(\tau-s)}P_c g_u R_a\ehts P_c v\|_{L^1}d\tau \nonumber \\
&+\int_s^{s+\frac{t-s}{4}} \|e^{-iH(t+s-2\tau)}P_c\|_{L^1\rightarrow L^\infty}\|\ehts P_c g_{\bar{u}}\overline{R}_a e^{iH(\tau-s)}  P_c \bar{v}\|_{L^1}d\tau \nonumber \\
&+\int_{s+\frac{t-s}{4}}^{t-\frac{t-s}{4}} \|e^{-iH(t-\tau)}P_c\|_{L^1\rightarrow L^infty}\|g_{\bar{u}} \overline{R}_a e^{iH(\tau-s)}  P_c \bar{v}\|_{L^1}d\tau \nonumber \\
&+\int_{t-\frac{t-s}{4}}^t \|e^{-iH(t+s-2\tau)}P_c\|_{L^1\rightarrow L^\infty}\|\ehts P_c g_{\bar{u}}\overline{R}_a e^{iH(\tau-s)}  P_c \bar{v}\|_{L^1}d\tau \nonumber \\
&\leq\int_s^t \frac{C}{|t-s|^{\frac{3}{2}}}\sup\|\widehat{g_u}\|_{L^1}\|v\|_{L^1}d\tau+\int_s^{s+\frac{t-s}{4}} \frac{C}{|t+s-2\tau|^{\frac{3}{2}}}\sup \|\widehat{g_{\bar{u}}}\|_{L^1}\|v\|_{L^1}d\tau \nonumber \\
&+\int_{s+\frac{t-s}{4}}^{t-\frac{t-s}{4}} \frac{C}{|t-\tau|^{\frac{3}{2}}}\|g_{\bar{u}}\|_{L^1}\|v\|_{L^1}d\tau+\int_{t-\frac{t-s}{4}}^t \frac{C}{|t+s-2\tau|^{\frac{3}{2}}}\sup \|\widehat{g_{\bar{u}}}\|_{L^1}\|v\|_{L^1}d\tau \nonumber \\ &\leq C\frac{\|v\|_{L^1}}{|t-s|^{\frac{1}{2}}}(\|g_{\bar{u}}\|_{L^1}+\sup(\|\widehat{g_u}\|_{L^1}+\|\widehat{g_{\bar{u}}}\|_{L^1})) \nonumber \end{align}
Now let us investigate the long time bevaviour of the forcing term $f(t)$. We will assume $t>s+1$ and seperate $f(t)$ into four parts as follows,

$$f(t)=\underbrace{\int_s^{s+\frac{1}{4}}\cdots}_{I_1}+\underbrace{\int_{s+\frac{1}{4}}^{t-\frac{1}{4}}\cdots}_{I_2}+\underbrace{\int_{t-\frac{1}{4}}^t\cdots}_{I_3}$$ We will start with $I_2$ for which we are away from the singularities around $\tau=s$ and $\tau=t$. Then for $I_1$ and $I_3$ we will use J-S-S type estimate to remove the singularities.
\begin{align}
\|I_2\|_{L^\infty}&\leq\int_{s+\frac{1}{4}}^{t-\frac{1}{4}} \|\eht P_c (g_u R_a\ehts P_c v+g_{\bar{u}}\overline{R}_a e^{iH(\tau-s)} P_c \bar{v})\|_{L^\infty}d\tau \nonumber \\
&\leq\int_{s+\frac{1}{4}}^{t-\frac{1}{4}} \frac{C}{|t-\tau|^{\frac{3}{2}}} (\|g_u R_a\ehts P_c v\|_{L^1}+\|g_{\bar{u}}\overline{R}_a e^{iH(\tau-s)} P_c \bar{v}\|_{L^1})d\tau \nonumber \\
&\leq\int_{s+\frac{1}{4}}^{t-\frac{1}{4}}\frac{C}{|t-\tau|^{\frac{3}{2}}}( \|g_u\|_{L^1}+\|g_{\bar{u}}\|_{L^1})\|\ehts P_c v\|_{L^\infty}d\tau \nonumber \\
&\leq\int_{s+\frac{1}{4}}^{t-\frac{1}{4}}\frac{C}{|t-\tau|^{\frac{3}{2}}}( \|g_u\|_{L^1}+\|g_{\bar{u}}\|_{L^1}) \frac{\|v\|_{L^1}}{|\tau-s|^{\frac{3}{2}}}d\tau \nonumber \\
&\leq C\frac{\|v\|_{L^1}}{|t-s|^{\frac{3}{2}}}(\|g_u\|_{L^1}+\|g_{\bar{u}}\|_{L^1}) \nonumber  \end{align} 

\begin{align}
\|I_1\|_{L^\infty}&\leq\int_s^{s+\frac{1}{4}} \|\ehs P_c e^{iH(\tau-s)} P_c g_u R_a\ehts P_c v\|_{L^\infty}d\tau \nonumber \\
&\quad+\int_s^{s+\frac{1}{4}} \|e^{-iH(t+s-2\tau)} \ehts P_c g_{\bar{u}}\overline{R}_a e^{iH(\tau-s)} P_c \bar{v}\|_{L^\infty}d\tau \nonumber \\
&\leq\int_s^{s+\frac{1}{4}}\frac{C}{|t-s|^{\frac{3}{2}}}\|e^{iH(\tau-s)} P_c g_u R_a\ehts P_c v\|_{L^1}d\tau \nonumber \\
&\quad+\int_s^{s+\frac{1}{4}}\frac{C}{|t+s-2\tau|^{\frac{3}{2}}}\|\ehts P_c g_{\bar{u}}\overline{R}_a e^{iH(\tau-s)} P_c \bar{v}\|_{L^1}d\tau \nonumber \\
&\leq\int_s^{s+\frac{1}{4}} C \Big(\frac{1}{|t-s|^{\frac{3}{2}}}+\frac{1}{|t-s-1|^{\frac{3}{2}}}\Big)\sup(\|\widehat{g_u}\|_{L^1}+\|\widehat{g_{\bar{u}}}\|_{L^1})\|v\|_{L^1}d\tau \nonumber \\
&\leq C\frac{\|v\|_{L^1}}{|t-s|^{\frac{3}{2}}} \nonumber \end{align}
\begin{align}
\|I_3\|_{L^\infty}&\leq\int_{t-\frac{1}{4}}^t \|\eht P_c g_u R_a e^{iH(t-\tau)} P_c e^{-iH(t-s)} P_cv\|_{L^\infty}d\tau \nonumber \\
&\quad+\int_{t-\frac{1}{4}}^t \|\eht P_c g_{\bar{u}}\overline{R}_a e^{iH(t-\tau)} P_c e^{iH(t+s-2\tau)} P_c \bar{v}\|_{L^\infty}d\tau \nonumber \\
&\leq C\frac{\|v\|_{L^1}}{|t-s|^{\frac{3}{2}}} \nonumber \end{align}

Now it remains to show that $L(s)W$ is bounded in $L^\infty$. Again to remove the singularities we will split the integral in different parts. Let us consider $s\leq t\leq s+1$,
\begin{align}
&L(s)W=\int_s^t \eht P_c F_1(R_a W(\tau)) d\tau\nonumber \\
&=\int_s^t \eht P_c g_u R_a\Big[\int_s^\tau e^{-iH(\tau-\tau')} P_c [F_1(R_a e^{-iH(\tau'-s)}P_c v)+F_1(R_a W(\tau'))]d\tau'\Big]d\tau \nonumber \\
&\quad+\int_s^t \eht P_c g_{\bar{u}}\overline{R}_a\Big[\int_s^\tau e^{iH(\tau-\tau')}P_c[\overline{ F_1(R_a e^{-iH(\tau'-s)}P_c v)}+\overline{ F_1(R_a W(\tau'))}]d\tau'\Big]d\tau \nonumber \end{align}
All the terms will be either of the following forms
$$L_1=\int_s^t \eht P_c g_u R_a\int_{s}^\tau e^{-iH(\tau-\tau')} P_c X(\tau')d\tau'd\tau$$
$$L_2=\int_s^t \eht P_c g_{\bar{u}}\overline{R}_a\int_{s}^\tau e^{iH(\tau-\tau')} P_c \overline{X(\tau')}d\tau'd\tau$$
where $X(\tau')=g_u R_a e^{-iH(\tau'-s)}P_c v, g_{\bar{u}}\overline{R}_a e^{iH(\tau'-s)}P_c \bar{v}, g_u R_a W(\tau'), g_{\bar{u}}\overline{R_a W(\tau')}$

In what follows we will add $e^{iH(t-\tau)}$ and $\eht$ terms after $g_u R_a$ and $g_{\bar{u}}\overline{R}_a$ then we will estimate the terms in a similiar way as we estimated $I_1$ and $I_3$.
\begin{equation}
L_1=\int_s^t \eht P_c g_u R_a e^{iH(t-\tau)}\int_{s}^\tau e^{-iH(t-\tau')} P_c X(\tau')d\tau'd\tau \label{l1} \end{equation}
\begin{equation}L_2=\int_s^t \eht P_c g_{\bar{u}}\overline{R}_a e^{iH(t-\tau)}\int_{s}^\tau e^{-iH(t-2\tau+\tau')} P_c \overline{X(\tau')}d\tau'd\tau \label{l2} \end{equation}
\begin{itemize}
\item For $X(\tau')=g_u R_a e^{-iH(\tau'-s)}P_c v$ we have
\begin{align}
\|L_1\|_{L^\infty}&\leq\int_s^t \|\eht P_c g_u R_a e^{iH(t-\tau)}\|_{L^\infty\rightarrow L^\infty}\int_s^{\tau}\|e^{-iH(t-s)}P_c\|_{L^1\rightarrow L^\infty} \|e^{iH(\tau'-s)}P_c g_u R_a e^{-iH(\tau'-s)}P_c v\|_{L^1}d\tau'd\tau \nonumber \\
&\leq\int_s^t \|\widehat{g_u}\|_{L^1}\int_s^\tau \frac{C}{|t-s|^{\frac{3}{2}}}\|\widehat{g_u}\|_{L^1}\|v\|_{L^1}d\tau'd\tau\leq C\sqrt{t-s}\|v\|_{L^1}\leq C\|v\|_{L^1}\quad\textrm{for }s\leq t\leq s+1 \nonumber \end{align}
\begin{align}
\|L_2\|_{L^\infty}&\leq\int_s^{t-\frac{t-s}{4}} \|\eht P_c\|_{L^1\rightarrow L^\infty}\|g_{\bar{u}}\overline{R}_a\int_s^{\tau}e^{-iH(\tau-s)}P_c e^{-iH(\tau'-s)}P_c g_{\bar{u}}\overline{R}_a e^{iH(\tau'-s)}P_c \bar{v}\|_{L^1}d\tau'd\tau \nonumber \\
&+\int_{t-\frac{t-s}{4}}^t \|\eht P_c g_{\bar{u}}\overline{R}_a e^{iH(t-\tau)}\|_{L^\infty\rightarrow L^\infty} \nonumber \\
&\quad\quad\quad\quad\times \int_s^{\tau}\|e^{-iH(t+s-2\tau)}P_c\|_{L^1\rightarrow L^\infty}\|e^{-iH(\tau'-s)}P_c g_{\bar{u}}\overline{R}_a e^{iH(\tau'-s)}P_c \bar{v}\|_{L^1}d\tau'd\tau \nonumber \\
&\leq\int_s^{t-\frac{t-s}{4}}\frac{C}{|t-\tau|^{\frac{3}{2}}} \|g_{\bar{u}}\|_{L^1}\int_s^\tau \frac{C}{|\tau-s|^{\frac{3}{2}}}\|\widehat{g_{\bar{u}}}\|_{L^1}\|v\|_{L^1}d\tau'd\tau \nonumber \\
&+\int_{t-\frac{t-s}{4}}^t \|\widehat{g_{\bar{u}}}\|_{L^1}\int_s^\tau \frac{C}{|t+s-2\tau|^{\frac{3}{2}}}\|\widehat{g_{\bar{u}}}\|_{L^1}\|v\|_{L^1}d\tau'd\tau\leq C\sqrt{|t-s|}\|v\|_{L^1} \nonumber \end{align}
\item For $X(\tau')=g_{\bar{u}}\overline{R}_a e^{iH(\tau'-s)}P_c \bar{v}$ we have
\begin{align}
\|L_1\|_{L^\infty}&\leq\int_s^t \|\eht P_c g_u R_a e^{iH(t-\tau)}\|_{L^\infty\rightarrow L^\infty}\nonumber \\ &\quad\quad\quad\times\Big[\int_{s}^{s+\frac{t-s}{4}} \|e^{-iH(t+s-2\tau')} P_c\|_{L^1\rightarrow L^\infty}\|e^{-iH(\tau'-s)}P_c g_{\bar{u}}\overline{R}_a e^{iH(\tau'-s)}P_c \bar{v}\|_{L^1}d\tau'\nonumber \\
&\quad\quad\quad+\int_{s+\frac{t-s}{4}}^{t-\frac{t-s}{4}}\|e^{-iH(t-\tau')}\|_{L^1\rightarrow L^\infty}\|g_{\bar{u}}\overline{R}_a e^{iH(\tau'-s)}P_c \bar{v}\|_{L^1}d\tau' \nonumber \\
&\quad\quad\quad+\int_{t-\frac{t-s}{4}}^{\tau}\|e^{-iH(t+s-2\tau')} P_c\|_{L^1\rightarrow L^\infty}\|e^{-iH(\tau'-s)}P_c g_{\bar{u}}\overline{R}_a e^{iH(\tau'-s)}P_c \bar{v}\|_{L^1}d\tau'\Big]d\tau \nonumber \\
&\leq\int_s^t \|\widehat{g_u}\|_{L^1}\Big[\int_s^{s+\frac{t-s}{4}}\frac{C}{|t+s-2\tau'|^{\frac{3}{2}}}\|\widehat{g_{\bar{u}}}\|_{L^1}d\tau'+\int_{s+\frac{t-s}{4}}^{t-\frac{t-s}{4}}\frac{C}{|t-\tau'|^{\frac{3}{2}}}\|g_{\bar{u}}\|_{L^1}\frac{C}{|\tau-s|^{\frac{3}{2}}}d\tau'\nonumber \\
&\quad\quad\quad\quad\quad\quad\quad\quad\quad+\int_{t-\frac{t-s}{4}}^\tau\frac{C}{|t+s-2\tau'|^{\frac{3}{2}}}\|\widehat{g_{\bar{u}}}\|_{L^1}d\tau'\Big]\|v\|_{L^1}d\tau\nonumber \\
&\leq C\sqrt{|t-s|}\|v\|_{L^1} \nonumber\end{align}
\begin{align}
\|L_2\|_{L^\infty}&\leq\int_s^t \|\eht P_c g_{\bar{u}}\overline{R}_a e^{iH(t-\tau)}\|_{L^\infty\rightarrow L^\infty}\Big[\int_{s}^{s+\frac{t-s}{4}}\| e^{-iH(t-2\tau+2\tau'-s)}e^{iH(\tau'-s)} P_c g_u R_a e^{-iH(\tau'-s)}P_c v\|_{L^\infty}d\tau' \nonumber \\
&\quad\quad\quad+\int_{s+\frac{t-s}{4}}^\tau \|e^{-iH(t-2\tau+\tau')} P_c\|_{L^1\rightarrow L^\infty} \|g_u R_a e^{-iH(\tau'-s)}P_c v\|_{L^1}d\tau' \Big]d\tau \nonumber \\
&\leq\int_s^t \|\widehat{g_{\bar{u}}}\|_{L^1}\Big[\int_s^{s+\frac{t-s}{4}} \frac{C}{|t-2\tau+2\tau'-s|^{\frac{3}{2}}}\|\widehat{g_{\bar{u}}}\|_{L^1} \|v\|_{L^1}d\tau'+\int_{s+\frac{t-s}{4}}^\tau \frac{C}{|t-2\tau+\tau'|^{\frac{3}{2}}}\|g_{\bar{u}}\|_{L^1} \frac{\|v\|_{L^1}}{|\tau'-s|^{\frac{3}{2}}}d\tau'\Big] d\tau \nonumber \\
&\leq C\sqrt{|t-s|}\|v\|_{L^1} \nonumber \end{align}

\item For $X(\tau')=g_u R_a W(\tau')$ and $g_{\bar{u}} \overline{R_a W(\tau')}$ we will change the order of integration,
\begin{align}
\|L_1\|_{L^\infty}&\leq\int_s^t\int_{\tau'}^t\|\eht P_c g_u R_a e^{iH(t-\tau)}\|_{L^\infty\rightarrow L^\infty}\|e^{-iH(t-\tau')}P_c\|_{L^1\rightarrow L^\infty} \|g_u R_a W(\tau')\|_{L^1}d\tau d\tau' \nonumber \\
&\leq\int_s^t\int_{\tau'}^t \|\widehat{g_u}\|_{L^1}\frac{C}{|t-\tau'|^{\frac{3}{2}}}\|\xs g_u\|_{L^2}\|W\|_{\Lsn}d\tau d\tau' \nonumber \\
&\leq\int_s^t\|\widehat{g_u}\|_{L^1}\frac{C}{|t-\tau'|^{\frac{1}{2}}}\|\xs g_u\|_{L^2}\frac{C\|v\|_{L^1}}{(1+|\tau'-s|)^{\frac{1}{2}}} d\tau' \nonumber \\
&\leq C\sqrt{|t-s|}\|v\|_{L^1} \nonumber \end{align}
\begin{align}
\|L_2\|_{L^\infty}&\leq\int_s^t\int_{\tau'}^{t-\frac{t-\tau'}{4}} \|\eht P_c\|_{L^1\rightarrow L^\infty}\|g_u R_a  e^{-iH(\tau-\tau')} P_c g_u R_a W(\tau')\|_{L^1}d\tau d\tau' \nonumber \\
&+\int_s^t\int_{t-\frac{t-\tau'}{4}}^t \|\eht P_c g_u R_a e^{iH(t-\tau)}\|_{L^\infty\rightarrow L^\infty}\| e^{-iH(t+\tau'-2\tau)} P_c g_u R_a W(\tau')\|_{L^\infty}d\tau d\tau' \nonumber \\
&\leq\int_s^t\int_{\tau'}^{t-\frac{t-\tau'}{4}}\frac{C}{|t-\tau|^{\frac{3}{2}}}\|g_{\bar{u}}\|_{L^2}\|e^{-iH(\tau-\tau')}P_c\|_{L^2\rightarrow L^2}\|g_u R_a W(\tau')\|_{L^2}d\tau d\tau' \nonumber \\
&+\int_s^t\int_{\tau'}^{t-\frac{t-\tau'}{4}} \|\widehat{g_u}\|_{L^1}\frac{C}{|t+\tau'-2\tau|^{\frac{3}{2}}}\|g_u R_a W(\tau')\|_{L^1}d\tau d\tau'  \nonumber \\
&\leq\int_s^t\frac{C}{|t-\tau'|^{\frac{1}{2}}}\|g_{\bar{u}}\|_{L^2}\|\xs g_u\|_{L^\infty}\|W(\tau')\|_{\Lsn}d\tau' \nonumber \\
&+\int_s^t \|\widehat{g_u}\|_{L^1}\frac{C}{|t-\tau'|^{\frac{1}{2}}}\|\xs g_u\|_{L^2}\|W(\tau')\|_{\Lsn}d\tau' \nonumber \\
&\leq C\sqrt{|t-s|}\|v\|_{L^1} \nonumber
\end{align}
\end{itemize}

Similarly we will investigate the long time behavior of the operator $L(s)$ for $t>s+1$.
$$L(s)W(t)=\underbrace{\int_s^{t-\frac{1}{4}}\cdots}_{L_3}+\underbrace{\int_{t-\frac{1}{4}}^t\cdots}_{L_4}$$
\begin{align}
\|L_3\|_{L^\infty}&\leq\int_s^{t-\frac{1}{4}} \|\eht P_c F_1(R_a W(\tau))\|_{L^\infty} d\tau \nonumber \\
&\leq\int_s^{t-\frac{1}{4}} \|\eht P_c\|_{L^1\rightarrow L^\infty}\| F_1(R_a W(\tau))\|_{L^1} d\tau \nonumber \\
&\leq\int_s^{t-\frac{1}{4}} \frac{C}{|t-s|^{\frac{3}{2}}}(\|\xs g_u\|_{L^2}+\|\xs g_{\bar{u}}\|_{L^2})\|W\|_{\Lsn} d\tau \nonumber \\
&\leq\frac{C}{|t-s|^{\frac{3}{2}}}\int_s^{t-\frac{1}{4}} \frac{\|v\|_{L^1}}{(1+|\tau-s|)^{\frac{3}{2}}} d\tau\leq\frac{C}{|t-s|^{\frac{3}{2}}}\|v\|_{L^1} \nonumber \end{align}
In $L_4$ we will plug in \eqref{sow} once more:
\begin{align}
L_4&=\int_{t-\frac{1}{4}}^t \eht P_c F_1(R_a W(\tau)) d\tau\nonumber \\
&=\int_{t-\frac{1}{4}}^t \eht P_c g_u R_a\Big[\int_s^\tau e^{-iH(\tau-\tau')} P_c [F_1(R_a e^{-iH(\tau'-s)}P_c v)+ F_1(R_a W(\tau'))]d\tau'\Big]d\tau \nonumber \\
&\quad+\int_{t-\frac{1}{4}}^t \|\eht P_c g_{\bar{u}}\overline{R}_a\Big[\int_s^{\tau} e^{iH(\tau-\tau')}P_c[\overline{F_1(R_a e^{-iH(\tau'-s)}P_c v)}+\overline{F_1(R_a W(\tau'))}]d\tau'\Big]d\tau \nonumber \end{align}
Again we will add $e^{iH(t-\tau)}$ and $\eht$ terms after $g_u R_a$ and $g_{\bar{u}}\overline{R}_a$. Then  all the terms will be similar to $L_1$, $L_2$, $(\ref{l1})-(\ref{l2})$ respectively. After seperating the the inside integrals into pieces, we will estimate short time step integrals exactly the same way we did short time behavior by using JSS estimate, and the other integrals will be estimated using the usual norms.
\begin{itemize}
 \item For $X(\tau')=g_u R_a e^{-iH(\tau'-s)}P_c v$ we have
\begin{align}
\|L_1\|_{L^\infty}&\leq\int_{t-\frac{1}{4}}^t \|\eht P_c g_u R_a e^{iH(t-\tau)}\|_{L^\infty\rightarrow L^\infty}\nonumber \\
&\quad\quad\Big[\int_s^{s+\frac{1}{4}}\|e^{-iH(t-s)}P_c\|_{L^1\rightarrow L^\infty} \|e^{iH(\tau'-s)}P_c g_u R_a e^{-iH(\tau'-s)}P_c v\|_{L^1}d\tau' \nonumber \\
&\quad\quad+\int_{s+\frac{1}{4}}^{t-\frac{1}{4}} \|e^{-iH(t-\tau')} P_c\|_{L^1\rightarrow L^\infty}\|g_u R_a e^{-iH(\tau'-s)}P_c v\|_{L^\infty}d\tau' \nonumber \\
&\quad\quad+\int_{t-\frac{1}{4}}^{\tau}\|e^{-iH(t-s)}P_c\|_{L^1\rightarrow L^\infty} \|e^{iH(\tau'-s)}P_c g_u R_a e^{-iH(\tau'-s)}P_c v\|_{L^1}d\tau'\Big]d\tau \nonumber \\
&\leq \frac{C\|v\|_{L^1}}{|t-s|^{\frac{3}{2}}} \nonumber \end{align}
\begin{align}
\|L_2\|_{L^\infty}&\leq\int_{t-\frac{1}{4}}^t \|\eht P_c g_{\bar{u}}\overline{R}_a e^{iH(t-\tau)}\|_{L^\infty\rightarrow L^\infty} \nonumber \\
&\quad\quad\quad\quad\times\Big[ \int_s^{s+\frac{1}{4}}\|e^{-iH(t+s-2\tau)}P_c\|_{L^1\rightarrow L^\infty}\|e^{-iH(\tau'-s)}P_c g_{\bar{u}}\overline{R}_a e^{iH(\tau'-s)}P_c \bar{v}\|_{L^1}d\tau' \nonumber \\
&\quad\quad\quad+\int_{s+\frac{1}{4}}^{t-\frac{1}{4}}\|e^{-iH(t-2\tau+\tau')}P_c\|_{L^1\rightarrow L^\infty}\|g_{\bar{u}} \overline{R}_a e^{iH(\tau'-s)}P_c \bar{v}\|_{L^1}d\tau' \nonumber \\
&\quad\quad\quad+\int_{t-\frac{1}{4}}^{\tau}\|e^{-iH(t+s-2\tau)}P_c\|_{L^1\rightarrow L^\infty}\|e^{-iH(\tau'-s)}P_c g_{\bar{u}}\overline{R}_a e^{iH(\tau'-s)}P_c \bar{v}\|_{L^1}d\tau'\Big]d\tau \nonumber \\
&\leq\frac{C\|v\|_{L^1}}{|t-s|^{\frac{3}{2}}} \nonumber \end{align}
\item For $X(\tau')=g_{\bar{u}}\overline{R}_a  e^{iH(\tau'-s)}P_c \bar{v}$ we have
\begin{align}
\|L_1\|_{L^\infty}\leq&\int_{t-\frac{1}{4}}^t \|\eht P_c g_u R_a e^{iH(t-\tau)}\|_{L^\infty\rightarrow L^\infty}\nonumber \\ &\times\Big[\int_{s}^{s+\frac{1}{4}} \|e^{-iH(t+s-2\tau')} P_c\|_{L^1\rightarrow L^\infty}\|e^{-iH(\tau'-s)}P_c g_{\bar{u}}\overline{R}_a e^{iH(\tau'-s)}P_c \bar{v}\|_{L^1}d\tau'\nonumber \\
&\quad\quad\quad+\int_{s+\frac{1}{4}}^{t-\frac{1}{4}}\|e^{-iH(t-\tau')}\|_{L^1\rightarrow L^\infty}\|g_{\bar{u}}\overline{R}_a e^{iH(\tau'-s)}P_c \bar{v}\|_{L^1}d\tau' \nonumber \\
&\quad\quad\quad+\int_{t-\frac{1}{4}}^{\tau}\|e^{-iH(t+s-2\tau')} P_c\|_{L^1\rightarrow L^\infty}\|e^{-iH(\tau'-s)}P_c g_{\bar{u}}\overline{R}_a e^{iH(\tau'-s)}P_c \bar{v}\|_{L^1}d\tau'\Big]d\tau \nonumber \\
\leq&\int_{t-\frac{1}{4}}^t \|\widehat{g_u}\|_{L^1}\Big[\int_s^{s+\frac{1}{4}}\frac{C}{|t+s-2\tau'|^{\frac{3}{2}}}\|\widehat{g_{\bar{u}}}\|_{L^1}d\tau'+\int_{s+\frac{1}{4}}^{t-\frac{1}{4}}\frac{C}{|t-\tau'|^{\frac{3}{2}}}\|g_{\bar{u}}\|_{L^1}\frac{C}{|\tau-s|^{\frac{3}{2}}}d\tau'\nonumber \\
&\quad\quad\quad\quad\quad\quad\quad\quad\quad+\int_{t-\frac{1}{4}}^\tau\frac{C}{|t+s-2\tau'|^{\frac{3}{2}}}\|\widehat{g_{\bar{u}}}\|_{L^1}d\tau'\Big]\|v\|_{L^1}d\tau\nonumber \\
\leq& C\frac{\|v\|_{L^1}}{|t-s|^{\frac{3}{2}}} \nonumber\end{align}
\begin{align}
\|L_2\|_{L^\infty}&\leq\int_s^t \|\eht P_c g_{\bar{u}}\overline{R}_a e^{iH(t-\tau)}\|_{L^\infty\rightarrow L^\infty}\Big[\int_{s}^{s+\frac{1}{4}}\| e^{-iH(t-2\tau+2\tau'-s)}e^{iH(\tau'-s)} P_c g_u R_a e^{-iH(\tau'-s)}P_c v\|_{L^\infty}d\tau' \nonumber \\
&\quad\quad\quad+\int_{s+\frac{1}{4}}^\tau \|e^{-iH(t-2\tau+\tau')} P_c\|_{L^1\rightarrow L^\infty} \|g_u R_a e^{-iH(\tau'-s)}P_c v\|_{L^1}d\tau' \Big]d\tau \nonumber \\
&\leq\int_s^t \|\widehat{g_{\bar{u}}}\|_{L^1}\Big[\int_s^{s+\frac{1}{4}} \frac{C}{|t-2\tau+2\tau'-s|^{\frac{3}{2}}}\|\widehat{g_{\bar{u}}}\|_{L^1} \|v\|_{L^1}d\tau'+\int_{s+\frac{1}{4}}^\tau \frac{C}{|t-2\tau+\tau'|^{\frac{3}{2}}}\|g_{\bar{u}}\|_{L^1} \frac{\|v\|_{L^1}}{|\tau'-s|^{\frac{3}{2}}}d\tau'\Big] d\tau \nonumber \\
&\leq \frac{C\|v\|_{L^1}}{|t-s|^{\frac{3}{2}}} \nonumber \end{align}
\item $L_1$, $L_2$ terms corresponding to $X(\tau')=g_u R_a W(\tau')$ and $g_{\bar{u}} \overline{R_a W(\tau')}$
\begin{align}
\|L_1\|_{L^\infty}&\leq\int_{t-\frac{1}{4}}^t \|\eht P_c g_u R_a e^{iH(t-\tau)}\|_{L^\infty\rightarrow L^\infty}\int_s^{\tau}\|e^{-iH(t-\tau')}P_c\|_{L^1\rightarrow L^\infty} \|g_u R_a W(\tau')\|_{L^1}d\tau'd\tau \nonumber \\
&\leq\int_{t-\frac{1}{4}}^t \|\widehat{g_u}\|_{L^1}\int_s^{\tau}\frac{C}{|t-\tau'|^{\frac{3}{2}}}\|\xs g_u\|_{L^2}\frac{C\|v\|_{L^1}}{(1+|\tau'-s|)^{\frac{3}{2}}}d\tau' d\tau \nonumber \\
&\leq\int_{t-\frac{1}{4}}^t\Bigg[\int_s^{\frac{t+s}{2}} \frac{C}{|t-\tau'|^{\frac{3}{2}}}\frac{C\|v\|_{L^1}}{(1+|\tau'-s|)^{\frac{3}{2}}}d\tau +\int_{\frac{t+s}{2}}^\tau \frac{C}{|t-\tau'|^{\frac{3}{2}}}\frac{C\|v\|_{L^1}}{(1+|\tau'-s|)^{\frac{3}{2}}}d\tau\Bigg] d\tau' \nonumber \\
&\leq \frac{C\|v\|_{L^1}}{|t-s|^{\frac{3}{2}}} \nonumber
\end{align}
\begin{align}
\|L_2\|_{L^\infty}&\leq\int_{t-\frac{1}{4}}^t \|\eht P_c g_{\bar{u}}\overline{R_a} e^{iH(t-\tau)}\|_{L^\infty\rightarrow L^\infty}\int_s^{\tau}\|e^{-iH(t+\tau'-2\tau)}P_c\|_{L^1\rightarrow L^\infty} \|g_u R_a W(\tau')\|_{L^1}d\tau'd\tau \nonumber \\
&\leq\int_{t-\frac{1}{4}}^t \|\widehat{g_{\bar{u}}}\|_{L^1}\int_s^{\tau}\frac{C}{|t+\tau'-2\tau|^{\frac{3}{2}}}\|\xs g_u\|_{L^2}\frac{C\|v\|_{L^1}}{(1+|\tau'-s|)^{\frac{3}{2}}}d\tau' d\tau \nonumber \\
&\leq\int_{t-\frac{1}{4}}^t\Bigg[\int_s^{\frac{\tau+s}{2}} \frac{C}{|t+\tau'-2\tau|^{\frac{3}{2}}}\frac{C\|v\|_{L^1}}{(1+|\tau'-s|)^{\frac{3}{2}}}d\tau +\int_{\frac{\tau+s}{2}}^\tau \frac{C}{|t-\tau'|^{\frac{3}{2}}}\frac{C\|v\|_{L^1}}{(1+|\tau'-s|)^{\frac{3}{2}}}d\tau\Bigg] d\tau' \nonumber \\
&\leq \frac{C\|v\|_{L^1}}{|t-s|^{\frac{3}{2}}} \nonumber
\end{align}
\end{itemize}

\noindent Now combining all the above estimates we get $$\|W(t)\|_{L^\infty}\leq\left\{ \begin{array}{ll}
C |t-s|^{\frac{1}{2}} & \textrm{for $|t-s|\leq 1$} \\
\frac{C}{|t-s|^{\frac{3}{2}}} & \textrm{for $|t-s|>1$} \end{array} \right.$$ This finishes the proof of $(ii)$.


$(iii)$  We split $f$ given by \eqref{sow}:
$$f=\underbrace{\int_s^{s+1}\cdots}_{I_1}+\underbrace{\int_{s+1}^t\cdots}_{I_2}$$
For $I_1$ integral, it suffices to show that $\|g_u(\tau)R_a\ehts P_c v\|_{L^2}\in L^1_\tau[s,s+1]$. Since $g_u(\tau)$ has bounded derivatives we have $\|g_u(\tau)-g_u(s)\|_{L^3}\leq C|\tau-s|$, then by H\"{o}lder inequality in space, $$\|(g_u(\tau)-g_u(s))R_a\ehts P_c v\|_{L^2}\leq\|g_u(\tau)-g_u(s)\|_{L^3}\|\ehts P_c v\|_{L^6} \in L^1_\tau$$ Now it suffices to show $\|g_u(s)R_a\ehts P_c v\|_{L^2}\in L^1_\tau$. For any $\tilde{v}\in L^2$ we have $$\|g_u(s)R_a\ehts P_c v\|_{L^2}=\<\tilde{v},g_u(s)R_a\ehts P_c v\>=\<e^{iH(\tau-s)}P_c g_u(s)R_a\tilde{v},v\>\leq\underbrace{\|e^{iH(\tau-s)}P_c R_a g_u(s)\tilde{v}\|_{L^6}}_{\in L^2_\tau}\|v\|_{L^{6/5}}$$ Since $L^2_\tau[s,s+1]\hookrightarrow L^1_\tau[s,s+1]$, $\|e^{iH(\tau-s)}P_c R_a g_u(s)\tilde{v}\|_{L^6}\in L^1_\tau$.
\begin{align}
\|I_2\|_{L^2}&\leq C\Big(\int_{s+1}^t \|g_u R_a\ehts P_c v\|_{L^{\rho'}}^{\gamma'}d\tau\Big)^{\frac{1}{\gamma'}}+C\Big(\int_{s+1}^t \|g_{\bar{u}}\overline{R}_a e^{iH(\tau-s)} P_c \bar{v}\|_{L^{\rho'}}^{\gamma'}d\tau\Big)^{\frac{1}{\gamma'}} \nonumber \\
&\leq C\Big(\int_{s+1}^t \|\xs g_u\|_{L^{\frac{3\gamma}{2}}}^{\gamma'}\|\ehts P_c v\|_{L_{-\sigma}^2}^{\gamma'}d\tau\Big)^{\frac{1}{\gamma'}}+C\Big(\int_{s+1}^t \|\xs g_{\bar{u}}\|_{L^{\frac{3\gamma}{2}}}^{\gamma'}\|e^{iH(\tau-s)} P_c \bar{v}\|_{L_{-\sigma}^2}^{\gamma'}d\tau\Big)^{\frac{1}{\gamma'}} \nonumber \\
&\leq C\Big(\int_{s+1}^t  \frac{d\tau}{|\tau-s|^{3(\frac{1}{2}-\frac{1}{p})\gamma'}}\Big)^{\frac{1}{\gamma'}}<\infty \nonumber
\end{align}
At the first inequality we used Strichartz estimate with $(\gamma,\rho)$ with $\gamma>2$ and the last inequality holds since $3(\frac{1}{2}-\frac{1}{p})\gamma'>1$ for $p=6$ and $\gamma>2$.
Similarly we will estimate $L(s)W$.
\begin{align}
\|L(s)W(t)\|_{L^2}&\leq C\Big(\int_s^t \|g_u R_a W+g_{\bar{u}}\overline{R_aW}\|_{L^{\rho'}}^{\gamma'}d\tau\Big)^{\frac{1}{\gamma'}} \nonumber \\
&\leq C\Big(\int_s^t \|\xs (g_u+g_{\bar{u}})\|_{L^{\frac{3\gamma}{2}}}^{\gamma'}\|W\|_{\Lsn}^{\gamma'}d\tau\Big)^{\frac{1}{\gamma'}}\nonumber \\
&\leq C\Big(\int_s^t \frac{1}{(1+|\tau-s|)^{3(\frac{1}{2}-\frac{1}{p})\gamma'}}\Big)^{\frac{1}{\gamma'}}<\infty \nonumber \end{align}
Hence $T(t,s):L^{p'}\rightarrow L^2$ is bounded for $p=6$. 
This finishes the proof of part $(iii)$ and the theorem. $\Box$


\section{Appendix}\label{se:ap}

\subsection{J-S-S type estimates}
In \cite{kn:JSS} the authors obtain the following
estimate\footnote{Their theorem is stated differently but the proof
can be easily adapted to obtain the advertised estimate}:
\begin{theorem}\label{th:jss} If $W:\mathbb{R}^n\mapsto\mathbb{C}$ has
Fourier transform $\widehat{W}\in L^1(\mathbb{R}^n)$ then for any
$t\in\mathbb{R}$ and any $1\le p\le\infty$ we have:
$$\|e^{-i\Delta t}We^{i\Delta t}\|_{L^p\mapsto L^p}\le \|\widehat{W}\|_{L^1}$$
\end{theorem}
In what follows we are going to generalize the estimate to the
semigroup of operators generated by $-\Delta+V:$
\begin{theorem}\label{th:gjss} Assume $V:\mathbb{R}^n\mapsto\mathbb{R}$
and $W:\mathbb{R}^n\mapsto\mathbb{C}$ have Fourier transforms in
$L^1(\mathbb{R}^n).$ Then for any $T>0$ there exist a constant $C_T$
independent of $W$ such that for any $-T\le t\le T$ and any $1\le
p\le\infty$ we have:
$$\|e^{-i(-\Delta +V) t}We^{i(-\Delta +V)t}\|_{L^p\mapsto L^p}\le C_T\|\widehat{W}\|_{L^1}.$$
One can choose $C_T=\exp(2\|\widehat{V}\|_{L^1}T).$
\end{theorem}
The proof relies on existence of finite time wave operators:
\begin{lemma}\label{lem:waveop}
If $V:\mathbb{R}^n\mapsto\mathbb{R}$ has Fourier transform in
$L^1(\mathbb{R}^n)$ then for any $T>0$ there exist a constant $C_T$
such that for any $-T\le t\le T$ and any $1\le p\le\infty$ we have:
$$\|e^{-i(-\Delta +V) t}e^{-i\Delta t}\|_{L^p\mapsto L^p}\le C_T\|\widehat{W}\|_{L^1},\quad
\|e^{i\Delta t}e^{i(-\Delta +V) t}\|_{L^p\mapsto L^p}\le
C_T\|\widehat{W}\|_{L^1}.$$ One can choose
$C_T=\exp(\|\widehat{V}\|_{L^1}T).$
\end{lemma}
{\bf Proof of Lemma:} Let
$$H=-\Delta+V$$ then $H$ is a self adjoint operator on $L^2$ with
domain $H^2,$ (note that $V\in L^\infty,$) hence it generates a
group of isometric operators:
$$e^{-iHt}:L^2\mapsto L^2,\qquad t\in\mathbb{R}.$$
Consequently:
\begin{equation}\label{def:qt}Q(t)=e^{-iHt}e^{-i\Delta
t}:L^2\mapsto L^2,\qquad t\in\mathbb{R},
\end{equation} is also a
family of isometric operators. Their infinitesimal generators are:
$$\frac{dQ}{dt}=-ie^{-iHt}Ve^{-i\Delta
t}=-i\underbrace{e^{-iHt}e^{-i\Delta t}}_{Q(t)}\
\underbrace{e^{i\Delta t}Ve^{-i\Delta t}}_{Q_0(t)}$$ Hence
\begin{equation}\label{eq:qt}
Q(t)=\mathbb{I}_d-i\int_0^tQ(s)Q_0(s)ds
\end{equation}
where $$Q_0(t)=e^{i\Delta t}Ve^{-i\Delta t}:L^p\mapsto L^p,\qquad
1\le p\le \infty$$ is bounded uniformly by $\|\widehat V\|_{L^1},$
see Theorem \ref{th:jss}.

The contraction principle shows that for any $T>0$ and any $1\le
p\le\infty$ the linear equation \eqref{eq:qt} has a unique solution
in the Banach space $C([-T,T],B(L^p,L^p)).$ Since on $L^2\bigcap
L^p$ the solution is given by \eqref{def:qt} and $L^2\bigcap L^p$ is
dense in $L^p$ we obtain that for any $-T\le t\le T$ and any $1\le
p\le\infty,$ $e^{-iHt}e^{-i\Delta t}$ has a unique extension to a
bounded operator on $L^p.$ Applying the $L^p$ norm in \eqref{eq:qt}
we get:
$$\|Q(t)\|_{L^p}\le 1+\int_0^t\|Q(s)\|_{L^p} \|\widehat V\|_{L^1}ds$$
and by Gronwall inequality:
$$\|Q(t)\|_{L^p}\le e^{\|\widehat V\|_{L^1} |t|}\le e^{\|\widehat
V\|_{L^1} T},\quad {\rm for} -T\le t\le T.$$

A similar argument can be made for $Q^*(t)=e^{i\Delta t}e^{iHt}.$

The Lemma is now completely proven.
\bigskip

\noindent {\bf Proof of Theorem \ref{th:gjss}:} For $H,$ $Q$ and
$Q^*$ as in the proof of the previous Lemma we have:
$$
e^{-iHt}WE^{iHt}=\underbrace{e^{-iHt}e^{-i\Delta
t}}_{Q(t)}\underbrace{e^{i\Delta t}We^{-i\Delta t}}_{L^p\mapsto L^p\
{\rm bounded}}\underbrace{e^{i\Delta t}e^{iHt}}_{Q^*(t)}.$$ Hence
using Theorem \ref{th:jss} and Lemma \ref{lem:waveop} we get for any
$1\le p\le\infty:$
$$\|e^{-iHt}WE^{iHt}\|_{L^p\mapsto L^p}\le e^{2\|\widehat
V\|_{L^1} T}\|\widehat W\|_{L^1},\qquad {\rm for}\ -T\le t\le T. $$

The theorem is now completely proven.

\begin{remark}\label{jss} To obtain the linear estimates in Section
\ref{se:lin} we used Theorem \ref{th:gjss} in the form:
$$\|e^{iHt}WR_ae^{-iHt}\|_{L^p\mapsto L^p}\le C\|\widehat
W\|_{L^1},\qquad {\rm for}\ 0\le t\le 1$$ where $W$ is the effective
potential induced by the nonlinearity, see next subsection, while
$R_a$ is the linear operator defined in Lemma \ref{le:pcinv}.
\end{remark}
To see why the above estimate holds consider $f\in L^p\bigcap
L^2\bigcap {\cal H}_0.$ Then by \eqref{radef} we have for a certain
$z=z(f)\in\mathbb{C}:$
$$e^{iHt}WR_ae^{-iHt}f=e^{iHt}We^{-iHt}f+ze^{iHt}W\psi_0.$$
Theorem \ref{th:gjss} applies directly to the first term on the
right hand side, while for the second term we use, see
\eqref{z:est}:
$$|z|\le 2\|f\|_{L^p}\sqrt{\left\|\frac{\partial \psi_E}{\partial
a_2}\right\|_{L^{p'}}+\left\| \frac{\partial \psi_E}{\partial
a_1}\right\|_{L^{p'}}},\qquad \frac{1}{p}+\frac{1}{p'}=1$$ and the
fact that $\psi_0$ is an e-vector of $H$ with e-value $E_0<0$ hence
$$\|e^{iHt}W\psi_0\|_{L^p}=\|e^{iHt}We^{-iHt}e^{iE_0t}\psi_0\|_{L^p}\le
C\|\widehat W\|_{L^1}\|\psi_0\|_{L^p},$$ where again we used Theorem \ref{th:gjss}.

\subsection{Smoothness of the effective potential}
In this section we will prove Proposition \ref{pr:dg} i.e. $\widehat{g'(\pe)}$ and $\widehat{(\frac{g(\pe)}{\pe})}$
\bigskip

\noindent From by Corollary \ref{co:decay}, we have $\pe\in H^2$ which implies $\pe\in L^p$ for $2\leq p\leq\infty$. Also from \eqref{gest}, by integrating, we get $|g'(s)|\leq C(|s|^{1+\p_1}+|s|^{1+\p_2})$. Hence $|g'(\pe)|\leq C(|\pe|^{1+\p_1}+|\pe|^{1+\p_2})\in L^2$ and $|g''(\pe)|\leq C(|\pe|^{\p_1}+|\pe|^{\p_2})\in L^\infty$. Now we have
\begin{align}
\|\widehat{g'(\pe)}\|_{L^1}&=\|\frac{1}{1+|\xi|^2}(1+|\xi|^2)\widehat{g'(\pe)}\|_{L^1}\nonumber \\
&\leq\|\frac{1}{1+|\xi|^2}\|_{L^2}\|(1+|\xi|^2)\widehat{g'(\pe)}\|_{L^2}\nonumber \\
&\leq C(\|\widehat{g'(\pe)}\|_{L^2}+\|\widehat{\Delta g'(\pe)}\|_{L^2}) \nonumber \\
&\leq C(\|\underbrace{g'(\pe)}_{\in L^2}\|_{L^2}+\|\Delta
g'(\pe)\|_{L^2}) \nonumber\end{align} So it suffices to show that
$\Delta g'(\pe)\in L^2$. Similarly it is enough to show that $\Delta
(\frac{g(\pe)}{\pe})\in L^2$.

\begin{equation}
\Delta g'(\pe)=g'''(\pe)|\nabla \pe|^2+\underbrace{g''(\pe)}_{\in
L^\infty}\underbrace{\Delta\pe}_{\in L^2}\label{eq:dgpe}
\end{equation} and
\begin{equation}
\Delta
(\frac{g(\pe)}{\pe})=(\frac{g''(\pe)}{\pe}-2\frac{g'(\pe)}{\pe^2}+2\frac{g(\pe)}{\pe^3})|\nabla
\pe|^2+(\underbrace{\frac{g'(\pe)}{\pe}-\frac{g(\pe)}{\pe^2}}_{\in
L^\infty})\underbrace{\Delta\pe}_{\in L^2}\label{eq:gpe}
\end{equation} We will use the following comparison theorem proved in
\cite[Theorem 2.1]{deiftsimon} to get the upper bound for the
$\nabla\pe$ and lower bound for $\pe$:
\begin{theorem}\label{th:comp} Let $\varphi\geq0$ be continuous on $\overline{\R^3\setminus K}$ and $A\geq B\geq0$ for some closed set $K$. Suppose that on $\overline{\R^3\setminus K}$, in the distributional sense, $$\Delta |\psi|\geq A|\psi|;\quad\quad \Delta \varphi\leq B\varphi$$ and that $|\psi|\leq\varphi$ on $\partial K$ and $\psi$, $\varphi\rightarrow0$ as $x\rightarrow\infty$. Then $|\psi|\leq\varphi$ on all of $\R^3\setminus K$.
\end{theorem}

\noindent Note that $\frac{\partial \pe}{\partial x_1}$ and $\pe$ are continuous and $\frac{\partial \pe}{\partial x_1},\ \pe\rightarrow0$ as $|x|\rightarrow\infty$. Hence $\left|\frac{g(\pe)}{\pe}\right|\leq C(|\pe|^{1+\p_1}+|\pe|^{1+\p_2})\rightarrow0$ as $x\rightarrow\infty$.
\bigskip

\noindent First we need the standard upper bound for $\pe\geq0$. For any $A<-E$, there exists $C_A$ depending on $A$ such that $\pe\leq C_A e^{-\sqrt{A}|x|}$. Indeed if $R$ is sufficiently large, on $\overline{\R^3\setminus B(0,R)}$ we have $$\Delta\pe=[-E+V(x)+\frac{g(\pe)}{\pe}]\pe\geq A\pe,\quad\text{and}\quad \Delta\varphi=A\varphi-\frac{2\sqrt{A}}{|x|}\varphi\leq A\varphi$$ and on $\partial B(0,R)$ we have $\pe\leq C_A e^{-\sqrt{A}|x|}$ for $C_A$ big enough. Then by Theorem \ref{th:comp} we have $\pe\leq C_A e^{-\sqrt{A}|x|}$ on $\R^3\setminus B(0,R)$.

\noindent To get the lower bound for $\pe$ we will choose $\varphi=\pe$ and $\psi=C e^{-\sqrt{A_2}|x|}$ in Theorem \ref{th:comp}. On $\overline{\R^3\setminus B(0,R)}$, fix $\varepsilon>0$, $A_2\geq-E+2\varepsilon$ and choose $R$ large enough such that  $\frac{2\sqrt{A_2}}{|x|}\leq\varepsilon$ for $|x|\geq R$. Then from \eqref{eq:ev} we have  $$\Delta \pe=[-E+V(x)]\pe+g(\pe)
\leq[-E+V+\frac{g(\pe)}{\pe}]\pe\leq (-E+\varepsilon)\pe$$ and for $A_2\geq -E+2\varepsilon$ we have  $$\Delta
\psi=A_2\psi-\frac{2\sqrt{A_2}}{|x|}\psi\geq(-E+\varepsilon)\psi$$
Choose $C$ such that $C
e^{-\sqrt{A_2}|x|}\leq\pe$ on $\partial B(0,R)$. Then by theorem
\ref{th:comp}, we have $C e^{-\sqrt{A_2}|x|}\leq\pe$ for $|x|>R$.

\noindent We will show that for $\psi=\frac{\partial \pe}{\partial x_1}$ and
$\varphi=C e^{-\sqrt{A_1}|x|}$ where $A_1<-E$ hypothesis of the
theorem \ref{th:comp} is satisfied.

Differentiating the eigenvalue equation \eqref{eq:ev} with respect to $x_1$ we
get $$\Delta \frac{\partial \pe}{\partial
x_1}=[-E+V(x)]\frac{\partial \pe}{\partial
x_1}+g'(\pe)\frac{\partial \pe}{\partial x_1}+\frac{\partial
V}{\partial x_1}\pe$$ Let $$f^{\pm}=\max\{ 0,\pm f\}\quad\text{and}\quad
S_{\leq}=\{x\in\R^3|\Big|\frac{\partial \pe}{\partial x_1}\Big|\leq\pe\}\quad\text{and}\quad
S_{\geq}=\{x\in\R^3|\Big|\frac{\partial \pe}{\partial x_1}\Big|\geq\pe\}$$ Fix $A_1<-E$, choose $R$ large enough such that $-E+V(x)+g'(\pe)-\Big|\frac{\partial V}{\partial x_1}\Big|\geq A_1$ on $|x|\geq R$. Let $S=S_{\leq}\cup B(0,R)$, then on $\overline{\R\setminus S}$ we have
$$\Delta \Big|\frac{\partial \pe}{\partial x_1}\Big|\geq A_1\Big|\frac{\partial \pe}{\partial x_1}\Big|$$ Now, by continuity of $\frac{\partial \pe}{\partial x_1}$ there exists $C_1$ such that $\Big|\frac{\partial \pe}{\partial x_1}\Big|e^{\sqrt{A_1}|x|}\leq C_1$ on $|x|=R$. Since both  on $\frac{\partial \pe}{\partial x_1}$ and $\pe$ are continuous we have $\Big|\frac{\partial \pe}{\partial x_1}\Big|=\pe\leq C_2 e^{\sqrt{A_1}|x|}$ on $\partial S_\leq$. So on $\partial S$, we have $\Big|\frac{\partial \pe}{\partial x_1}\Big|\leq\max\{C_1,C_2\} e^{\sqrt{A_1}|x|}$

\noindent Therefore by
theorem \ref{th:comp}, we have $|\nabla \pe|\leq C
e^{-\sqrt{A_1}|x|}$

Now we can prove Proposition \ref{pr:dg}

\textbf{Proof of Proposition \ref{pr:dg}} By {\bf (H2')} we have $|g'''(s)|<\frac{C}{s^{1-\alpha_1}}+Cs^{\alpha_2 -1},$ $s>0,$  $0<\alpha_1\le\alpha_2;$ then $$|g'''(\pe)|\nabla
\pe|^2|\leq\frac{C}{\pe^{1-\p_1}}|\nabla \pe|^2$$ and
$$|(\frac{g''(\pe)}{\pe}-2\frac{g'(\pe)}{\pe^2}+2\frac{g(\pe)}{\pe^3})|\nabla
\pe|^2|\leq\frac{C}{\pe^{1-\p_1}}|\nabla \pe|^2$$

Using the estimates for $|\nabla \pe|$ and $\pe$ and choosing
$2\sqrt{A_1}>\sqrt{A_2}$, we get that $\Delta g'(\pe),\Delta
(\frac{g(\pe)}{\pe})\in L^2$. Hence we get the desired estimates for
$\widehat{g'(\pe)}$ and $\widehat{\frac{g(\pe)}{\pe}}$.

\bibliographystyle{plain}
\bibliography{ref}
\end{document}